\documentclass{article}

\usepackage{amsmath,amscd}

\usepackage{epsfig}

\usepackage{amssymb}
\usepackage[usenames,dvipsnames]{color}

\usepackage{latexsym}

\include{plaquettemacro}

\usepackage{hyperref}
 
\newtheorem{theorem}{Theorem}[section]
\newtheorem{prop}{Proposition}[section]
 \newtheorem{lemma}[theorem]{Lemma}

\numberwithin{equation}{section}

\newcommand{\mbf}{{\mathbb  F}}

\newcommand{\mbr}{{\mathbb R}}

\newcommand{\mpm}{{\mathcal P}M}

\newcommand{\mpx}{{\mathcal P}X}

\newcommand{\oab}{\omega_{(A,B)}}

\newcommand{\ova}{{\overline{a}}}

\newcommand{\ovA}{{\overline{A}}}

\newcommand{\ovAb}{  {\overline{\mathbf A}}             }

\newcommand{{\tlg}}{\tilde\gamma }

\newcommand{{\tlG}}{\tilde\Gamma }

\newcommand{\catmbt}{{\mathbb P}^{\rm bt}(M)}

\newcommand{\catmpt}{{\mathbb P}^{\rm bt}(P)}

\newcommand{\pap}{{\mathcal P}_{\ovA}P}

\newcommand{\mbbc}{{\mathbf C}}
\newcommand{\mbbz}{\mathbf {Z} }

\newcommand{\mbba}{\mathbf {A }}
\newcommand{\mbbb}{\mathbf {B} }
\newcommand{\mbg}{\mathbf {G} }

\newcommand{\mbk}{\mathbf {K} }

\newcommand{\mbv}{\mathbf {V} }

\newcommand{\mbp}{\mathbf {P} }
\newcommand{\Obj}{{\rm Obj }}
\newcommand{\Mor}{{\rm Mor }}

\begin{document}

\title{Pathspace Connections and Categorical Geometry  }


\author{Saikat Chatterjee\\
Department of Mathematics \\
 Harish Chandra Research Institute\\
Chhatnag Road,Jhusi,Allahabad-211 019,\\  
Uttar Pradesh, INDIA \\
email:saikat.chat01@gmail.com\\
\and Amitabha Lahiri\\ S.~N.~Bose National Centre for Basic Sciences \\ Block JD,
  Sector III, Salt Lake, Kolkata 700098 \\
  West Bengal, INDIA \\
 email: amitabhalahiri@gmail.com
 \and
  Ambar N. Sengupta\\  Department of Mathematics\\
  Louisiana State University\\  Baton
Rouge, Louisiana 70803, USA\\
email:   ambarnsg@gmail.com 
 }

\maketitle

\begin{abstract}  We study a type of connection forms, given by Chen integrals, over pathspaces by placing such forms within a category-theoretic framework of principal bundles and connections.    We introduce a notion of `decorated' principal bundles, develop parallel transport on such bundles, and explore specific examples in the context of pathspaces.  \end{abstract}



\section{Introduction}\label{intro}

In this paper we develop a theory of categorical geometry and explore   specific examples involving geometry over spaces of paths.  Our first objective is to develop a framework that encodes special properties, such as parametrization-invariance, of connection forms on path spaces. For this paper we focus on the case of connections over pathspaces given by first-order Chen integrals.  We then develop a theory of `decorated' principal bundles and parallel-transport in such bundles. These constructions all sit naturally inside a framework of categorical connections that we develop. A background motivation is to develop a framework that provides a unified setting for both ordinary
gauge theory, governing interactions between point particles, and higher gauge theory, governing the interaction of string-like, or higher-dimensional, extended objects.

We now summarize our results through an overview of the paper.  All categories we work with are `small': the objects and morphisms form sets. 

\begin{itemize}
\renewcommand{\labelitemi}{$\bullet$}
\item  In section \ref{S:connpath} we study a connection form $\oab$ on a principal bundle over a pathspace.   This connection form is invariant under reparametrization of the paths. We then describe the connection form $\oab$:
\begin{equation}\label{E:defoab2intro}
\omega_{(A,B)}({\tilde v})= A\bigl({\tilde v}(t_0)\bigr)+\tau\left[ \int_{t_0}^{t_1}B\left({\tilde \gamma}'(t), {\tilde v}(t)\right)\,dt\right],
\end{equation}
specifying the meaning of all the terms involved here. Briefly and roughly put, $\oab$ is a $1$-form on pathspace obtained by a `point-evaluation' of a traditional $1$-form (the first term on the right) and a first-order Chen integral (the second term on the right)). Next in Proposition \ref{P:omegareparinv} we prove that this connection form is also invariant under reparametrization of paths, and so induces a connection form on the space of reparametrization-equivalence classes of paths. In Proposition \ref{P:omegareparinv} we show that $\oab$ does have properties analogous to traditional connection forms on bundles.

\item  In section \ref{S:ptcs} we consider equivalence classes of paths, identifying paths that differ from each other by erasure of backtracked segments (that is, a composite path $c_2{\overline a}ac_1$, where $\overline a$ is the reverse of $a$, is considered equivalent to $c_2c_1$).  The results of section \ref{S:ptcs}, especially Theorem 
\ref{T:CLSpathT2.1ext}, show that $\oab$ specifies a connection form on the space of backtrack-erasure equivalence classes of paths. This is the result that links the geometry  with category theory: it makes it possible to view $\oab$ as specifying a connection form over a category whose morphisms are backtrack-erased paths. In Theorem \ref{T:thinhomot} we prove that $\oab$ respects another common way of identifying paths: paths $\gamma_1$ and $\gamma_2$ are said to be `thin homotopic' if one can obtain $\gamma_2$ from $\gamma_1$ by means of a homotopy that `wiggles' $\gamma_1$  along itself. Thus, $\oab$ specifies a connection also over the space of thin-homotopy equivalence classes of paths.

\item    In section \ref{S:2grp} we study  the notion of a  {\em categorical group}; briefly, this is a category whose object set and morphism set are both groups. The main result, Theorem \ref{T:catgrp2grp} establishes the equivalence between categorical groups and {\em crossed modules} specified by pairs of groups $(G,H)$. These results are  known in the literature but we feel it is useful to present this coherent account, as there are many different conventions and definitions used in the literature and our presentation in this section provides us with notation, conventions, and results for use in later sections.  We also include several examples here.

\item Section \ref{S:categorbundle} introduces the key notion of a {\em principal categorical bundle}  ${\mathbf P}\to {\mathbf B}$, with a categorical group ${\mathbf G}$ as `structure group ' and with both ${\mathbf P}$ and ${\mathbf B}$ being categories.  This general framework does not require ${\mathbf B}$ to have a smooth structure. We give examples, including one that uses backtrack-erased pathspaces. We conclude the section by showing that the notion of `reduction' of a principal bundle carries over to this categorical setting. (In this work we do not explore categorical analogs of local triviality, a topic that is central to most other works in this area.)

\item   In section \ref{S:pcatpath}  we introduce the notion of a {\em decorated} principal bundle. This gives a  useful example of a categorical principal bundle whose structure depends on the action of a given  crossed module. Briefly, we start with an ordinary principal bundle  $\pi:P\to B$ equipped with a connection  $\ovA$; we form a categorical principal bundle whose   base category  ${\mathbf B}$ has object set $B$  and backtrack-erased paths on $B$ as morphisms; the bundle category ${\mathbf P}$ has object set $P$ and morphisms of the form $({\tilde\gamma},h)$, with ${\tilde\gamma}$ being any $\ovA$-horizontal path on $P$ and $h$ running over a group $H$. Thus the morphisms are horizontal paths `decorated' with elements of $H$. The result is a categorical principal bundle whose structure categorical group is specified by the pair of groups $G$ and $H$, with $G$ acting on objects of ${\mathbf P}$ and a semi-direct product of $G$ and $H$ acting on the morphisms of ${\mathbf P}$.

\item   In section \ref{S:catcon} we introduce the notion of a {\em categorical connection} on a categorical principal bundle. We present several examples, and then show, in Theorem \ref{T:connectpapbt},  how to construct a categorical connection on the bundle of decorated paths and then, in Theorem \ref{T:absdecconn},  in a more abstractly decorated categorical bundle.

\item  In section \ref{S:pathcat} we describe,  in a precise way, categories whose morphisms are (equivalence classes of) paths, noting that there are two natural composition laws for morphisms.

\item In section \ref{S:catconpath} we construct categorical connections at a `higher' geometric level: here the objects are paths, and the morphisms are paths of paths. 

\item We present our final and most comprehensive example of a categorical connection in section\ref{S:ptdp}, where we develop parallel-transport of `decorated' paths over a space of paths. Thus the transport of a decorated path $({\tilde\gamma}, h)$ is specified through the data $({\tilde\Gamma}, h,k)$, where $\tilde\Gamma$ is a path of paths on the bundle, horizontal with respect to a pathspace connection $\oab$,  $h\in H$ decorates the initial  (or source)  path ${\tilde\gamma}_0=s({\tilde\Gamma})$ for ${\tilde\Gamma}$, and $k\in K$ encodes the rule for producing the resulting final decorated path $({\tilde\gamma}_1, h_1)$.
\item Section \ref{S:asocpt} presents a brief account of associated bundles, along with parallel-transport in such bundles, in the categorical framework. 

 \end{itemize}

 There is a considerable and active literature at the juncture of category theory and geometry. Without attempting to review the existing literature we make some remarks.  The relationship between higher categorical structures and quantum theories was first explored extensively by Freed \cite{Freed}. Works broadly related to ours include those of Baez et al. \cite{BaezDol, BaezLauda, BaezSchr, BaezSchr2, BaezWise}, Bartels \cite{Bart},   Schreiber and Waldorf \cite{SW1,SW2}, Abbaspour and Wagemann \cite{AbbasWag},  and Viennot \cite{Vien}. We study neither local trivialization of $2$-bundles nor the related notion of gerbes, both of which are explored in the other works just cited. Our definition of categorical connection in section \ref{S:catcon} appears to be new; other notions, such as $2$-connections \cite{BaezSchr2} may be found in the literature. Roughly speaking, our approach stays much closer to geometry than category theory in comparison to much of the $2$-geometry literature; this is because our primary objective is to create a framework for the  specific type of connections  given by (\ref{E:defoab2intro}).

 Before proceeding, let us briefly review some language connected with the standard framework of connections on principal bundles for ease of reference when considering the categorical definitions we introduce later. 
 
   For a Lie group $G$, a principal $G$-bundle $\pi:P\to B$  is a smooth surjection, where $P$ and $B$ are smooth manifolds, along with a free right action of $G$ on $P$:
$$P\times G\to P:(p,g)\mapsto R_g(p)=pg,$$
with $\pi(pg)=\pi(p)$ for all $p\in P$ and $g\in G$,
and there is local triviality: every point of $B$ has a neighborhood $U$ and there is a diffeomorphism
$$\phi:U\times G\to \pi^{-1}(U)$$
such that $\phi(u,g)h=\phi(u,gh)$  and $\pi\phi (u,g)=u$ for all $u\in U$ and $g,h\in G$.  In the categorical formulation we will develop in section \ref{S:categorbundle} we will not (in the present paper) impose local triviality.

A vector $v\in T_pP$ is {\em vertical} if $d\pi_p(v)=0$. 
A connection $\omega$ is a $1$-form on $P$ with values in the Lie algebra $L(G)$, and
satisfies the following conditions: (i)  $\omega(R_g'(p)X)={\rm Ad}(g^{-1})\omega(X)$ for all $p\in P$, $X\in T_p$ and $g\in G$, and (ii) $\omega({\tilde Y}(p))=Y$ for all $p\in P$ and $Y\in L(G)$, where ${\tilde Y}(p)=\frac{d}{dt}|_{t=0}pe^{tY}$.  The crucial consequence of this definition is that  $\omega$ decomposes each $T_pP$ as a direct sum of the vertical subspace $\ker d\pi_p$ and the {\em horizontal subspace} $\ker\omega_p$ in a manner `consistent' with the action of $G$. This leads to the notion of {\em parallel transport}: if $\gamma:[t_0,t_1]\to B$ is a   $C^1$ path and $u_0$ is a point on the fiber $\pi^{-1}\bigl(\gamma(t_0)\bigr)$ then there is a unique {\em horizontal lift} ${\tilde \gamma}_{\omega}:[t_0,t_1]\to P$, that is a $C^1$ path for which
\begin{equation}\label{E:omegacprim0}
\omega\left({\tilde \gamma}_{\omega}'(t)\right)=0\qquad\hbox{for all $t\in [t_0,t_1]$,}\end{equation}
that initiates at $u_0$.  It is in terms of the notion of parallel transport that we state the definition of a {  categorical connection } on a categorical bundle in section \ref{S:catcon}.  The {\em curvature} $2$-form $F^\omega$, containing information on infinitesimal holonomies, is defined to be
\begin{equation}\label{E:defFomega}
F^\omega=d\omega +\frac{1}{2}[\omega,\omega]\end{equation}
where, for any $L(G)$-valued $1$-forms $\eta$ and $\zeta$, the $2$-form $[\eta,\zeta]$ is given by
$$[\eta,\zeta](X,Y)=[\eta(X),\zeta(Y)]-[\eta(Y),\zeta(X)]$$
for all  $X,Y\in T_pP$ and all $p\in P$.

\section{A connection form over a pathspace}\label{S:connpath}

In this section we prove certain invariance properties of certain connection forms over spaces of paths. Our results are proved for piecewise $C^1$ paths, but for the sake of consistency with results in later sections we frame the  discussions (as opposed to the formal statements and proofs)   in terms of $C^\infty$ paths.

   For a manifold $X$  the set  of all parametrized paths on $X$ is
\begin{equation} \bigcup_{t_0,t_1\in\mbr, t_0<t_1}C^\infty\left([t_0,t_1]; X\right),\end{equation}
   and we denote by $\mpx$ the quotient set obtained by identifying paths that differ by a constant time-translation reparametrization. Thus, $\gamma:[t_0,t_1]\to X$ is identified with $\gamma_{+a}:[t_0-a,t_1-a]\to X:t\mapsto \gamma(t+a)$ in $\mpx$.  We will often not make a notational distinction between $\gamma$ and its equivalence class $[\gamma]$ of such time-translation reparametrizations, and   in using the term `parametrized path' we will   often not distinguish between time-translation reparametrizations of the same path.

We   work with a connection $\ovA$ on a principal $G$-bundle 
$$\pi:P\to M,$$
where $G$ is a Lie group.  We denote by
\begin{equation}\label{E:defpap}
\pap,\end{equation}
 the set of all $\ovA$-horizontal $C^\infty$ paths $\gamma:[t_0,t_1]\to P$,  again identifying paths that differ by a constant time-translation. (In   \cite{CLSpath} we used this notation but without  any identification procedure.) For a given connection $\ovA$ we view $\pap$ in a natural  but informal way as a principal $G$-bundle over $\mpm$; the projection $\pap\to\mpm$ is the natural one induced by $\pi:P\to M$, and the right action of $G$ on $\pap$ is also given simply from the action of $G$ on $P$.

Consider a path of $\ovA$-horizontal paths, specified through a $C^\infty$map
$${\tilde\Gamma}:[s_0,s_1]\times [t_0,t_1]\to P:(s,t)\mapsto {\tilde\Gamma}(s,t)={\tilde\Gamma}_s(t) ,$$
  where $s_0<s_1$ and $t_0<t_1$,  such that ${\tilde\Gamma}_s$ is $\ovA$-horizontal for every $s\in [s_0,s_1]$. In \cite[2.1]{CLSpath} it was shown that 
\begin{equation}\label{def:vTPAPinteg} 
 \ovA({\tilde
v}(t))= \ovA({\tilde
v}(t_0)) +\int_{t_0}^t F^{\ovA}\left({\tilde\gamma}'(s),
{\tilde v}(s)\right)\,ds \qquad \hbox{for all $t\in [t_0,t_1]$.}  
\end{equation}
where ${\tilde v}:[t_0,t_1]\to TP: t\mapsto \partial_s{\tilde\Gamma}(s_0,t)$ is the vector field along ${\tilde\gamma}={\tilde\Gamma}_{s_0}$ pointing in the variational direction. We refer to (\ref{def:vTPAPinteg}) as a  {\em tangency condition}.

Note that (\ref{def:vTPAPinteg}) is meaningful even if the path ${\tilde\gamma}$ is piecewise $C^1$ and the vector field ${\tilde v}$ merely continuous. Condition (\ref{def:vTPAPinteg}) is equivalent to   the requirement that
\begin{equation}\label{def:vTPAP} 
\begin{frac}{\partial\ovA({\tilde
v}(t))}{\partial t}\end{frac}= F^{\ovA}\left({\tilde\gamma}'(t),
{\tilde v}(t)\right)   
\end{equation}
hold at every point where ${\tilde\gamma}'(t)$ exists (which is all $t\in [t_0,t_1]$ except for finitely many points).

In the following result we show that the tangency condition (\ref{def:vTPAPinteg}) is independent of parametrizations, in the sense that  if ${\tilde v}$ and ${\tilde\gamma}$  satisfy  (\ref{def:vTPAPinteg}) then any reparametrization of ${\tilde\gamma}$ along with  the corresponding reparametrization of ${\tilde v}$ also satisfy (\ref{def:vTPAPinteg}). We denote  by
\begin{equation}\label{E:tangpap}
T_{\tilde\gamma}\pap.
\end{equation}
the vector space of all  vector fields  ${\tilde v}$ along ${\tilde\gamma}$ satisfying (\ref{def:vTPAPinteg}), with time-translates being identified.

 \begin{prop}\label{P:tangvecpathspace}  Let  $\ovA$ be a connection on a principal $G$-bundle, and let ${\tilde v} $ be a  continuous vector field along  an   $\ovA$-horizontal piecewise $C^1$ path ${\tilde\gamma}:[t_0,t_1]\to P$ satisfying (\ref{def:vTPAPinteg}).  Then:
 \begin{itemize}
 \item[\rm (i)] for any strictly increasing   piecewise $C^1 $   bijective reparametrization function $\phi:[t'_0,t'_1]\to [t_0,t_1]$, the vector field  ${\tilde v}\circ\phi$, along the path ${\tilde\gamma}\circ\phi$   satisfies (\ref{def:vTPAPinteg}) with the domain $[t_0,t_1]$ replaced by $[t'_0,t'_1]$;
 \item[\rm (ii)] if $\pi\circ{\tilde \gamma}$ is constant on an open subinterval of $[t_0,t_1]$ then ${\tilde\gamma}$ is constant on that same subinterval, and if the vector field $v=\pi_*{\tilde v}$ is constant on that subinterval  then so is ${\tilde v}$.
 \end{itemize}
 \end{prop}
 
 Part (ii) stresses how strongly the behavior of the path $\pi\circ{\tilde\gamma}$ controls the behavior of a vector field $\tilde v$ that satisfies the tangency condition 
 (\ref{def:vTPAPinteg}).

 \noindent\underline{Proof}. (i) Since $\phi$ is a strictly increasing bijection $[t'_0,t'_1]\to[t_0,t_1]$ we have, in particular, $\phi(t'_0)=t_0$. From  (\ref{def:vTPAPinteg}) we have:
 \begin{equation}\begin{split}
 \ovA\left({\tilde v}\circ\phi(t')\right) &= \ovA\left({\tilde v}\left(t_0\right)\right) +\int_{t_0}^{\phi(t')}F^{\ovA}\left({\tilde\gamma}'(u),
{\tilde v} (u)\right)\,du\\
&=\ovA\left({\tilde v}(t_0)\right) + \int_{ t'_0 }^{  t' }F^{\ovA}\left({\tilde\gamma}'\bigl(\phi(s)\bigr),
{\tilde v} \bigl(\phi(s)\bigr)\right)\,\phi'(s)ds\quad\hbox{(setting $u=\phi(s)$)}\\
&= \ovA\left({\tilde v}\bigl(\phi(t'_0)\bigr)\right) + \int_{t'_0}^{t'}F^{\ovA}\left(({\tilde\gamma}\circ\phi)'(s),
{\tilde v}\circ \phi(s)\right)\, ds.
 \end{split}\end{equation}
 In the second line above the change of variables $u=\phi(s)$ can be done interval by interval on each of which $\phi$ has a continuous derivative.
 
 (ii) Now suppose $\gamma=\pi\circ{\tilde\gamma}$ is constant on an open interval $(T,T+\delta)\subset [t_0,t_1]$. Then, for any $t\in (T,T+\delta)$  the `horizontal   projection' $\pi_*{\tilde\gamma}'(t)$  is $0$ and the `vertical part' $\ovA\bigl({\tilde\gamma}'(t)\bigr)$ is also $0$ because ${\tilde\gamma}$ is $\ovA$-horizontal. Hence ${\tilde\gamma}'(t)=0$ for all such $t$ and hence ${\tilde\gamma}$ is constant on this interval.  Next suppose $v=\pi_*{\tilde v}$ is constant over $(T,T+\delta)$. Then from (\ref{def:vTPAPinteg}) we have
 \begin{equation}
 \ovA\bigl({\tilde v}(t)\bigr)-  \ovA\bigl({\tilde v}(T)\bigr)= \int_{T}^t F^{\ovA}\left({\tilde\gamma}'(s),
{\tilde v}(s)\right)\,ds =0
 \end{equation}
 for all $t\in (T,T+\delta)$ because of the constancy of ${\tilde\gamma}$. Thus, 
 $$\ovA\bigl({\tilde v}(t)\bigr)=  \ovA\bigl({\tilde v}(T)\bigr)$$
 for all $t\in (T,T+\delta)$. Thus, since the horizontal parts of ${\tilde v}(t)$ and ${\tilde v}(T)$ are also assumed to be equal (by constancy of $v$ over $(T,T+\delta)$) it follows that ${\tilde v}(t)$ equals ${\tilde v}(T)$ for all $t\in (T, T+\delta)$. \fbox{QED}

We will study connections on the bundle $\pap$, and on related bundles, using  crossed mdules: a  {\em {crossed module}}  $(G,H, \alpha, \tau)$ consist of groups $G$ and $H$, along with homomorphisms
\begin{equation}\begin{split}
\tau:H\to G & \qquad \hbox{and}\qquad \alpha: G\to {\rm Aut}(H).\end{split}\end{equation}
satisfying the  identities:
\begin{equation}\label{E:Peiffer}
\begin{split}
\tau\bigl(\alpha(g)h\bigr) &= g   \tau(h)  g^{-1}\\
\alpha\bigl(\tau(h)\bigr)(h') &=hh'h^{-1}
\end{split}
\end{equation}
for all $g\in G$ and $h\in H$.    

The structure of a {crossed module} was  introduced in 1949 by J. H. C. Whitehead \cite[sec 2.9]{Whit}, foreshadowed in the dissertation of Peiffer \cite{Peif} (sometimes misspelled as `Pfeiffer' in the literature).

When working in the category of Lie groups, we require that $\tau$ be smooth, and the map 
\begin{equation}\label{E:alpgaghh}
G\times H\to H:(g,h)\mapsto \alpha(g,h)=\alpha(g)h
\end{equation}
be smooth.  In this case we call $(G, H,\alpha,\tau)$ a {\em Lie crossed module}.

In what follows we work with a Lie crossed module  $(G, H, \alpha,\tau)$, connections forms $A$ and $\ovA$    on a principal $G$-bundle $\pi:P\to M$. We denote the right action of $G$ on the bundle space $P$ by
$$P\times G\to P:(p,g)\mapsto pg=R_gp$$
and the corresponding derivative map by
$$T_pP\to T_{pg}P:v\mapsto vg\stackrel{\rm def}{=}dR_g|_p(v)$$
for all $p\in P$ and $g\in G$. We also work with
 an $L(H)$-valued $2$-form $B$ on $P$ which is $\alpha$-equivariant, in the sense that
\begin{equation}\label{E:Balphequiv}
B_{pg}(vg,wg)=\alpha(g^{-1})B_p(v,w)\qquad\hbox{for all $g\in G$, $p\in P$, and $v,w\in T_pP$,}
\end{equation}
and vanishes on vertical vectors in the sense that
\begin{equation}\label{E:Bvert}
B(v,w)=0\quad\hbox{whenever  $v,w\in T_p$ and $\pi_*v=0$,}
\end{equation}  
for any $p\in P$.

For a piecewise $C^1$ $\ovA$-horizontal path ${\tilde\gamma}:[t_0,t_1]\to P$ and any continuous vector field ${\tilde v}$ along ${\tilde\gamma}$ we define
\begin{equation}\label{E:defoab2}
\omega_{(A,B)}({\tilde v})= A\bigl({\tilde v}(t_0)\bigr)+\tau\left[ \int_{t_0}^{t_1}B\left({\tilde \gamma}'(t), {\tilde v}(t)\right)\,dt\right].
\end{equation}
(In \cite{CLSpath} we used ${\tilde v}(t_1)$ instead of ${\tilde v}(t_0)$  in the definition of $\oab$, which is equivalent to changing $B$ to $B-F^\ovA$ here.) We  view $\oab$ as a $1$-form on a space of paths on $P$, expressed as 
 \begin{equation}\label{E:defomeAB}
\omega_{(A,B)}={\rm ev}_{t_0}^*A+\tau(Z)
\end{equation}
where ${\rm ev}_t$ is the evaluation map
$${\rm ev}_t:\pap\to P:  {\tilde \gamma}\mapsto{\tilde \gamma}(t)$$
and $Z$ is the $L(H)$-valued $1$-form on $\pap$ given by the Chen integral
\begin{equation}\label{E:defZ}
Z=\int_{\tilde\gamma}B.
\end{equation}

(The preceding definitions, but with ${\rm ev}_{t_0}$ replaced by  ${\rm ev}_{t_1}$, were introduced  in our earlier work \cite{CLSpath}.)

The first observation for $\oab$ is that it remains invariant when paths are reparametrized;  thus, it can be viewed as a $1$-form defined on the space of equivalence classes of paths, with the equivalence relation being reparametrization:

\begin{prop}\label{P:omegareparinv}  Let  $(G, H, \alpha,\tau)$ be a  Lie {crossed module}, and $\ovA$ be  a connection on a principal $G$-bundle $\pi:P\to M$.  Let ${\tilde\gamma}:[t_0,t_1]\to P$ be a   piecewise $C^1$ $\ovA$-horizontal path in $P$, and ${\tilde v}:[t_0,t_1]\to P$ be a continuous vector field along ${\tilde\gamma}$ that satisfies the constraint (\ref{def:vTPAPinteg}). Let $A$ be a connection form on $\pi:P\to M$, $B$ an $L(H)$-valued $\alpha$-equivariant $2$-form on $P$ vanishing on vertical vectors, and $\oab$ be as in (\ref{E:defoab2}). Then 
\begin{equation}\label{E:omegaabreparinv}
\omega_{(A,B)}({\tilde v})=\omega_{(A,B)}({\tilde v}\circ\phi)
\end{equation}
for any     strictly increasing  piecewise  $C^1$
bijective function $\phi:[t'_0,t'_1]\to [t_0,t_1]$, where on the left ${\tilde v}$ is along ${\tilde\gamma}$  and on the right ${\tilde v}\circ\phi$ is along ${\tilde\gamma}\circ\phi$. 

\end{prop} 

The significance of this result is that  it implies that $\omega_{(A,B)}$ can be viewed as an $L(G)$-valued $1$-form on the bundle space $\pap$; in fact, though, for convenience,  we defined $\pap$ using only constant-time translations, $\oab$ would be defined on a space of paths, with paths that are reparametrizations of each other by increasing piecewise $C^1$ bijective functions.

\noindent\underline{Proof}. The  proof is by  change-of-variables  along the lines of the proof of Proposition \ref{P:tangvecpathspace} (i):
\begin{equation}\begin{split}
\omega_{(A,B)}({\tilde v}\circ\phi) &= A\bigl(({\tilde v}\circ \phi)(t'_0)\bigr)\bigr)+\tau\left[ \int_{t'_0}^{t'_1}B\left(\bigl({\tilde \gamma}\circ\phi\bigr)'(u), ({\tilde v}\circ \phi)(u)\right)\,du\right]\\
&= A\bigl({\tilde v}(t_0)\bigr)+ \tau\left[ \int_{\phi(t_0)}^{\phi(t_1)}B\left( \phi'(u) {\tilde\gamma}'\bigl(\phi(u)\bigr), {\tilde v}\bigl(\phi(u)\bigr)
\right)\,du\right]\\
&=A\bigl({\tilde v}(t_0)\bigr)+ \tau\left[ \int_{\phi(t_0)}^{\phi(t_1)}B\left({\tilde\gamma}'\bigl(\phi(u)\bigr), {\tilde v}\bigl(\phi(u)\bigr)
\right)\,\phi'(u) du\right]\\
&=A\bigl({\tilde v}(t_0)\bigr)+ \tau\left[ \int_{t_0}^{t_1}B\left( {\tilde\gamma}'(t), {\tilde v}\bigr)
\right)\,dt\right]\quad\hbox{(switching to $t=\phi(u)$)}\\
&=\omega_{(A,B)}({\tilde v}).   \end{split}\end{equation}
\fbox{QED}

Next we check that $\omega_{(A,B)}$ has the properties of a connection form on a principal $G$-bundle:

\begin{prop}\label{P:omegaconnect}  Let  $(G, H, \alpha,\tau)$ be a  Lie {crossed module},   $\ovA$ and $A$ be   connection forms on a principal $G$-bundle $\pi:P\to M$.  Let $B$ be an $L(H)$-valued $\alpha$-equivariant $2$-form on $P$ vanishing on vertical vectors, and $\oab$ be as in (\ref{E:defoab2}).  Then
\begin{equation}\label{E:oABvg}
\oab({\tilde v}g)={\rm Ad}(g^{-1})\oab({\tilde v})
\end{equation}
for all $g\in G$ and all continuous vector fields ${\tilde v}$ along any  $\ovA$-horizontal, piecewise $C^1$,  path ${\tilde\gamma}:[t_0,t_1]\to P$. Moreover, if $Y$ is any element of the Lie algebra $L(G)$ and ${\tilde Y}$ is the vector field along ${\tilde\gamma} $ given by ${\tilde Y}(t)=\frac{d}{du}\big|_{u=0}{\tilde\gamma}(t)\exp(uY)$, then
\begin{equation}\label{E:oABH}
\oab({\tilde Y})=Y.
\end{equation}
\end{prop} 
\noindent\underline{Proof}. Both of these statements follow directly from the defining relation
\begin{equation}\label{E:defoab3}
\omega_{(A,B)}({\tilde v})= A\bigl({\tilde v}(t_0)\bigr)+\tau\left[ \int_{t_0}^{t_1}B\left({\tilde \gamma}'(t), {\tilde v}(t)\right)\,dt\right].
\end{equation}
If in this equation we use ${\tilde v}g$ instead of ${\tilde v}$ then the first term on the right is conjugated by ${\rm Ad}(g^{-1})$ because $A$ is a connection form, whereas in the second term the integrand is conjugated by $\alpha(g^{-1})$ (keeping in mind that ${\tilde v}g$ is   vector field along ${\tilde\gamma}g$), and then using the relation $\tau\bigl(\alpha(g^{-1})h\bigr) =g^{-1}\tau(h)g$ for all $h\in H$,  from which we have 
$$\tau\bigl(\alpha(g^{-1})Z\bigr) ={\rm Ad}(g^{-1})Z\qquad\hbox{for all $Z\in L(H)$.}$$
This shows  that the second term is also conjugated by ${\rm Ad}(g^{-1})$. Next, if ${\tilde v}={\tilde Y}$, where $Y\in L(G)$ then the  integrand in the second term on the right in (\ref{E:defoab3}) is $0$ because $B$ vanishes on vertical vectors and the first term is $Y$ (since $A$ is a connection form). \fbox{QED}

The existence and uniqueness of horizontal lifts relative to $\oab$ follow by using standard results on the existence and uniqueness of solutions of ordinary differential equations in Lie groups. The key observation needed here is that  to obtain the $\oab$-horizontal lift ${\tilde\Gamma}$ of a given path of paths $\Gamma: [s_0,s_1]\times[t_0,t_1]\to M:(s,t)\mapsto \Gamma_s(t)$, and a given initial path ${\tilde\Gamma}_{s_0}$,  one need only specify the motion $s\mapsto {\tilde\Gamma}_s(t_0)\in P$ and then each $\ovA$-horizontal path $t\mapsto {\tilde\Gamma}_s(t)$ is completely specified through the  initial value ${\tilde\Gamma}_s(t_0)\in P$.

\section{Parallel transport and backtrack equivalence}\label{S:ptcs}

In this section we prove that   pathspace parallel-transport by the connection $\omega_{(A,B)}$ induces, in a natural way, parallel-transport over the space of backtrack-erased paths, a notion that we will explain drawing on ideas from the careful   treatment by L\'evy \cite{Le10}.

We say that a map $\gamma:[t_0,t_1]\to X$ into a space $X$  {\em backtracks over $[T, T+\delta]$}, where $t_0\leq T<T+2\delta<t_1$, if 
\begin{equation}\label{E:gammaTdubt}
\gamma(T+u)=\gamma(T+2\delta-u)\qquad\hbox{for all $u\in [0,\delta]$.}
\end{equation}
By {\em erasure of the backtrack} over $[T, T+\delta]$ from $\gamma$ we will mean the map
\begin{equation}\label{E:gammbterases}
 [t_0,t_1-2\delta]: t\mapsto \begin{cases} \gamma(t) & \hbox{if $t\in [t_0,T]$;}\\
\gamma(t-2\delta) & \hbox{if $t\in [T+2\delta, t_1]$.}
\end{cases}
\end{equation}
Clearly, this is continuous if $\gamma$ is continuous, and is piecewise $C^1$ if $\gamma$ is piecewise $C^1$.
(In later sections we will work with $C^\infty$ paths that are constant near the initial and final points.  When working with the class of such paths  we consider only backtrack erasures that preserve the $C^\infty$ nature.)

 In the following we identify a parametrized path $c_1$ with a parametrized path $c_2$ if  $c_2=c_1\circ\phi$, where $\phi$ is a strictly increasing piecewise $C^1$ mapping of the domain of $c_2$ onto the domain of $c_1$. We say that  piecewise $C^1$  paths  $c_1$ and $c_2$ on $M$ are  {\em elementary backtrack  equivalent}
  if there are   piecewise $C^1$  parametrized paths $a, b, d$ on $M$, and a strictly increasing piecewise $C^1$ function $\phi$ from some closed interval onto the domain of $c_1$,  and a strictly increasing piecewise $C^1$ function $\psi$ from some closed interval onto the domain of $c_2$, such that 
  \begin{equation}\label{E:abcd}
\{c_1\circ\phi , c_2\circ\psi\}=\{add^{-1}b, ab\},\end{equation}
where juxtaposition of paths denotes composition of paths, and for any path $d:[s_0,s_1]\to M$ the reverse $d^{-1}$ denotes the path
\begin{equation}\label{E:drev}d^{-1}:[s_0,s_1]\to M: s\mapsto d\bigl(s_1-(s-s_0)\bigr).\end{equation}
Condition (\ref{E:abcd}) means that one of the $c_i$ is obtained from the other by erasing the backtracking part $dd^{-1}$.

We say that  piecewise $C^1$ paths $b$ and $c$  are {\em backtrack equivalent}, denoted
\begin{equation}\label{E:defbtcb}
c\simeq_{\rm bt}b,\end{equation}
 if   there is a sequence of   paths $c=c_0, c_1,\ldots, c_n=b$, where $c_i$ is elementarily backtrack equivalent to $c_{i+1}$  for $i\in\{0,\ldots, n-1\}$.   
  
Backtrack equivalence is   an equivalence relation, and   if two paths are  backtrack equivalent then so are any of their reparametrizations. 
 Furthermore,  if
  $$a\simeq_{\rm bt} c\quad\hbox{and}\quad b\simeq_{\rm bt} d$$
and if the composite $ab$ is defined then so is $cd$ and
\begin{equation}\label{E:abapbp}
ab\simeq_{\rm bt} cd.\end{equation} 

By a {\em backtrack-erased path} $\gamma$ we will mean the  backtrack equivalence class  containing the specific path $\gamma$.

A tangent vector to a path $\gamma:[t_0,t_1]\to M$    is normally taken to be a vector field $v$, of suitable degree of smoothness, that has the property that $v(t)\in T_{\gamma(t)}M$, that is, it is a vector field {\em along} $\gamma$. Since we wish to identify backtrack equivalent paths we should not allow the vector field $v$ to `pry open' the backtracked parts of $\gamma$; thus we define a tangent vector $v$ to   $\gamma$  in the space of backtrack-equivalence classes of paths, to be a continuous (or suitably smooth) vector field $v$ along $\gamma:[t_0,t_1]\to M$, constant near $t_0$ and $t_1$, with the property that if $\gamma$   backtracks over $[T, T+\delta]$,  then $v$ also has the same backtrack:
\begin{equation}\label{E:btvts}v(T+u)= v(T+2\delta - u)\qquad\hbox{for all $u\in [0,\delta]$.}\end{equation}

\begin{prop}\label{P:backtracktangent} Let  $\ovA$ be a connection form on a principal $G$-bundle $\pi:P\to M$, where $G$ is a Lie group. 
Consider an $\ovA$-horizontal piecewise $C^1$ path ${\tilde\gamma}:[t_0,t_1]\to P$ and  a continuous vector field ${\tilde v}:[t_0,t_1]\to TP$    along $\tilde\gamma$  that  satisfies the tangency condition (\ref{def:vTPAP}). Suppose that the path $\gamma=\pi\circ\tilde\gamma$     backtracks over $[T, T+\delta]$ and suppose
the vector field $v=\pi_*{\tilde v}$ along $\gamma$ also  backtracks over  $[T, T+\delta]$. Then:
\begin{itemize}
\item[{\rm (i)}]  The   path $  \tilde\gamma$   backtracks over $[T, T+\delta]$.
\item[{\rm (ii)}]  The vector field ${\tilde v}$ backtracks over $[T, T+\delta]$.
\item[{\rm (iii)}]  Erasing the backtrack over $[T,T+\delta]$ from ${\tilde v}$  produces a vector field along the path obtained by erasing the backtrack over $[T,T+\delta]$ from ${\tilde\gamma}$ that continuous to satisfy the tangency condition (\ref{def:vTPAP}).
\end{itemize}
\end{prop}

In the following it is useful to keep in mind that a vector ${\tilde w}\in T_pP$ is completely determined by its  projection $\pi_*{\tilde w}\in T_{\pi(p)}M$ and its `vertical' part $\ovA\bigl({\tilde w}\bigr)\in L(G)$ (which can be transfered to an actual vertical vector in $T_pP$ by means of the right action of $G$ on $P$).

\noindent\underline{Proof}.  (i)  The backtracking in $\tilde\gamma$ follows because parallel-transport along the reverse of  a path is the reverse of the parallel-transport along the path.  Hence
$${\tilde\gamma}(T+u)={\tilde\gamma}(T+2\delta-u)\qquad\hbox{for all $u\in [0, \delta]$.}$$
Note that as consequence the velocity on the way back is minus the velocity on the way out:
\begin{equation}\label{E:gammabackvel}
{\tilde\gamma}'(T+u)=-{\tilde\gamma}'(T+2\delta-u)\qquad\hbox{for all $u\in [0, \delta]$.}\end{equation}

(ii) The vector ${\tilde v}(t)$ is  completely determining by its projection $v(t)=\pi_*{\tilde v}(t)$ and the `vertical part' $\ovA({\tilde
v}(t))\in L(G)$. To  prove that  ${\tilde v}$ backtracks over $[T,T+\delta]$ we need therefore only show that $\ova\circ {\tilde v}$ backtracks over $[T,T+\delta]$.

Recall the  tangency condition  (\ref{def:vTPAP}) :
\begin{equation}\label{E:tangcond}
\ovA({\tilde
v}(t))= \ovA({\tilde
v}(t_0)) +\int_{t_0}^t F^{\ovA}\left({\tilde\gamma}'(s),
{\tilde v}(s)\right)\,ds \qquad \hbox{for all $t\in [t_0,t_1]$.}  \end{equation}
 Since $F^\ovA$ vanishes on vertical vectors, on the right we can replace ${\tilde v}(t)$ by ${\tilde v}_h(t)$, the $\ovA$-horizontal lift of $v(t)$ as a vector in $T_{{\tilde\gamma}(t)}P$.  Hence
 
\begin{equation}\label{E:ovAv}
\ovA\bigl({\tilde v}(t)\bigr)=\ovA\bigl({\tilde v}(t_0)\bigr)+\int_{t_0}^tF^{\ovA}\big({\tilde\gamma}'(s), {\tilde v}_h(s)\bigr)\,ds.
\end{equation}
Clearly, the horizontal lift  ${\tilde v}_h $ backtracks to reflect the backtracking (\ref{E:btvts}) of $v$:
$${\tilde v }_h(T+u)={\tilde v}_h(T+2\delta-u)\qquad\hbox{for all $u\in [0,\delta]$.}$$
This, together with (\ref{E:gammabackvel}), implies that in the integral on the right in 
(\ref{E:ovAv}), the contribution from $s=T+u$ cancels the contribution from $s=T+2\delta-u$, for all  $u\in [0, \delta]$. Hence 
$$\ovA\bigl({\tilde v}(T+u)\bigr)=\ovA\bigl({\tilde v}(T+2\delta-u)\bigr)\qquad\hbox{for all $u\in [0, \delta]$.}$$

(iii) When $t\geq T+2\delta$ the  part of the integral on the right in (\ref{E:ovAv}) over $[T, T+2\delta]$ is completely erased, and so  $\ovA\bigl({\tilde v}(t)\bigr)$ retains  the same value if the backtrack were erased from $\gamma$. From this it follows that the tangency condition (\ref{E:ovAv}) continuous to hold when the backtracks over $[T,T+\delta]$ is erased from ${\tilde v}$ and from ${\tilde\gamma}$.
\fbox{QED}

We resume working with a {crossed module} $(G, H,\alpha,\tau)$, connections $A$, $\ovA$, and an $L(H)$-valued $2$-form $B$ on a principal $G$-bundle $\pi:P\to M$,   satisfying $\alpha$-equivariance (\ref{E:Balphequiv}) and vanishing on verticals (\ref{E:Bvert}). 
The following is a remarkable feature of the connection form $\omega_{(A,B)}$, showing that it specifies parallel-transport over the space of backtrack equivalence classes of paths:

\begin{theorem}\label{T:CLSpathT2.1ext} Let $(G,H, \alpha,\tau)$ be a Lie {crossed module},  $\ovA$ and $A$ be connections on  a principal $G$-bundle  $\pi:P\to M$, and let $\oab$ be the $1$-form given by  (\ref{E:defomeAB}).    
Then $\oab$ is well-defined as a $1$-form on tangent vectors to backtrack equivalence classes of $\ovA$-horizontal paths on $P$ in the following sense. Suppose $\gamma:[t_0,t_1]\to M$ is a piecewise $C^1$ path, and $\gamma_0$ a path obtained by erasing a backtrack over $[T,T+\delta]$ from $\gamma$.  Let $\tilde\gamma$ be an $\ovA$-horizontal lift of $\gamma$ and ${\tilde\gamma}_0$ the $\ovA$-horizontal lift of $\gamma_0$ with the same initial point as ${\tilde\gamma}$.  Suppose  ${\tilde v}$ is a continuous vector field along  $\tilde\gamma$  that backtracks over $[T, T+\delta]$. Then 
\begin{equation}\label{E:oabvv0}
\oab({\tilde v})=\oab({\tilde v}_0),
\end{equation}
where ${\tilde v}_0$ is the vector field   along ${\tilde\gamma}_0$ induced by $\tilde v$, through restriction. In particular, $\oab({\tilde v})
$ is $0$ if and only if  $\oab({\tilde v}_0)$ is also $0$.
\end{theorem}

\noindent\underline{Proof}.  Consider an $\ovA$-horizontal path ${\tilde\gamma}:[t_0,t_1]\to P$.  Suppose $\gamma=\pi\circ{\tilde\gamma}$  backtracks over $[T, T+\delta]$.  By Proposition \ref{P:backtracktangent}(i), 
  $\tilde\gamma$ also backtracks  over $[T, T+\delta]$.  Hence
\begin{equation}
{\tilde\gamma}'(T+u)  =-{\tilde\gamma}'(T+2\delta-u)\end{equation}
for all $u\in [0,\delta]$ for which either side exists. From the expression
\begin{equation}\label{E:defomeAB22}
\omega_{(A,B)}({\tilde v})= A\left({\tilde v}(t_0)\right)+\tau\left(\int_{t_0}^{t_1}B\left({\tilde \gamma}'(t), {\tilde v}(t)\right)\,dt
\right)
\end{equation}
we see then that in the second term on the right in the integral $\int_{t_0}^{t_1}\cdot\,dt$ we have a cancellation:
$$\int_{T}^{T+\delta} B\left({\tilde \gamma}'(t), {\tilde v}(t)\right)\,dt+
\int_{T+\delta}^{T+2\delta} B\left({\tilde \gamma}'(t), {\tilde v}(t)\right)\,dt=0.$$
Note also that the endpoint of a path is unaffected by backtrack erasure. Hence the first term on the right in (\ref{E:defomeAB22}) remains unchanged if  the backtrack over $[T, T+\delta]$ in $\tilde\gamma$ is erased. Thus the value $\omega_{(A,B)}({\tilde v})$ remains unchanged if the backtrack of $\gamma$ over any subinterval of $[t_0,t_1]$ is erased. In particular, ${\tilde v}$ is $\omega_{(A,B)}$-horizontal if and only if the backtrack-erased version of ${\tilde v}$ is $\omega_{(A,B)}$-horizontal.
\fbox{QED}

 There is a notion, known as `thin homotopy' equivalence, that is broader than backtrack equivalence. To keep technicalities from overpowering the ideas, we focus on the case of smooth paths and deformations.  Roughly speaking a thin homotopy wriggles a path back and forth along itself with no `transverse' motion.  The following result says essentially that parallel-transport  by $\oab$ is trivial along a thin homotopy. 
 
 \begin{theorem}\label{T:thinhomot} Let $(G,H, \alpha,\tau)$ be a Lie {crossed module},   $A$ and $\ovA$ connection forms on a principal  $G$-bundle $\pi:P\to M$,  and let $\oab$ be as given by  (\ref{E:defomeAB}).  Consider a $C^\infty$ map
 $$\Gamma:[s_0,s_1]\times[t_0,t_1]\to M:(u,v)\mapsto\Gamma_u(v)$$
for which $\partial_u\Gamma $ and $\partial_v\Gamma$ are linearly dependent at each point of $[s_0,s_1]\times[t_0,t_1]$, and, furthermore, $\Gamma$ keeps each endpoint of each $u$-fixed line $\{u\}\times[t_0,t_1]$ constant: 
\begin{equation}\label{E:Gammasti}
\Gamma(u,t_i)=\Gamma(s_0,t_i) \quad\hbox{ for all $u\in [s_0,s_1]$ and $i\in\{0,1\}$. }
\end{equation}
 Next let  ${\tilde\Gamma}_0:[t_0,t_1]\to P$ be any  $\ovA$-horizontal path that projects to the initial path $\Gamma_0$ on $M$. Then the $\oab$-horizontal lift of $u\mapsto\Gamma_u$, with ${\tilde\Gamma}_{s_0}={\tilde\Gamma}_0$, is the path
\begin{equation}\label{E:uGammathin}
[s_0,s_1]\to\pap:u\mapsto {\tilde\Gamma}_u\end{equation}
where ${\tilde\Gamma}_u:[t_0,t_1]\to P$ is the $\ovA$-horizontal lift of $\Gamma_u$
with initial point
\begin{equation}\label{E:Gamut0thin}
{\tilde\Gamma}_u(t_0)={\tilde\Gamma}_0(t_0) \end{equation}
for all $u\in [s_0,s_1]$.  Moreover, $\partial_u{\tilde\Gamma}$ and $\partial_v{\tilde\Gamma} $ are linearly dependent at each point of $[s_0,s_1]\times[t_0,t_1]$, where
${\tilde\Gamma}$ is the map $[s_0,s_1]\times[t_0,t_1]\to P:(u,v)\mapsto {\tilde\Gamma}_u(v)$.
 
 \end{theorem} 
 
 \noindent\underline{Proof}.    From the tangency condition (\ref{def:vTPAPinteg}) we have
 \begin{equation}\label{E:ovaGamthin} 
 \begin{split}
 \ovA(\partial_u{\tilde
\Gamma}(u,w)) &= \ovA(\partial_u{\tilde
\Gamma}(u,t_0)) \\
&\hskip .25in +\int_{t_0}^w F^{\ovA}\left(\partial_v{\tilde\Gamma}(u,v), \partial_u{\tilde\Gamma}(u,v)\right)\,dv \quad \hbox{for all $w\in [t_0,t_1]$.} 
\end{split} 
\end{equation}
The first term on the right   is $0$ because ${\tilde\Gamma}(\cdot, t_0)$ is constant by assumption (\ref{E:Gamut0thin}).  We now show that the second term is also zero.   By assumption $\partial_u\Gamma$ and $\partial_v\Gamma$ are linearly  dependent; thus there exist $a(u,v), b(u,v)\in\mbr$, not both zero, such that
\begin{equation}\label{E:abGammthin}
a(u,v)\partial_v{ \Gamma}(u,v) + b(u,v)\partial_u{ \Gamma}(u,v)=0,\end{equation}
for all $(u,v)\in [s_0,s_1]\times[t_0,t_1]$. Observe that
 $$ \pi_*\partial_v{\tilde\Gamma} =\partial_v(\pi\circ{\tilde\Gamma})=\partial_v\Gamma$$
 and, similarly,
 $$\pi_* \partial_u{\tilde\Gamma} =\partial_u\Gamma.$$
 Hence
\begin{equation}\label{E:abuvGamthin}
a(u,v)\partial_v{\tilde\Gamma}(u,v) + b(u,v)\partial_u{\tilde\Gamma}(u,v)\in\ker\pi_*.\end{equation}
(We will see shortly that the left side is in fact zero.) 
 Thus, up to addition of a vertical vector (on which $F^\ovA$ vanishes), the vectors $\partial_v{\tilde\Gamma}(u,v) $ and $\partial_u{\tilde\Gamma}(u,v) $ are linearly dependent and so the $2$-form $F^\ovA$ is $0$ when evaluated on this pair of vectors. 
 Hence
 \begin{equation}\label{E:ovaGam0thin} 
 \ovA(\partial_u{\tilde
\Gamma}(u,v))=0,
\end{equation}
for all $(u,v)\in [a_0,a_1]\times[t_0,t_1]$.  Using this and the $\ovA$-horizontality of ${\tilde\Gamma}_u$ we have
\begin{equation}\label{E:abuvtGamthin}
\ovA\left(a(u,v)\partial_v{\tilde\Gamma}(u,v) + b(u,v)\partial_u{\tilde\Gamma}(u,v)\right)=0,
\end{equation}
which, in combination with (\ref{E:abuvGamthin}), implies that
\begin{equation}\label{E:abuvtGam0thin}
 a(u,v)\partial_v{\tilde\Gamma}(u,v) + b(u,v)\partial_u{\tilde\Gamma}(u,v) =0.
\end{equation}
This proves that $\partial_u{\tilde\Gamma}$ and $\partial_v{\tilde\Gamma}$ are linearly dependent.

  From the definition of $\oab$ in (\ref{E:defoab2}) we have
 
  \begin{equation}\label{E:oabthin1}
 \oab(\partial_u{\tilde\Gamma})=A\bigl({\partial_u\tilde\Gamma}(u,t_0)\bigr)+\tau\left[\int_{t_0}^{t_1} B\left(\partial_v{\tilde\Gamma}(u,v), \partial_u{\tilde\Gamma}(u,v)\right)\,dv\right].
 \end{equation}
 The  second term vanishes because of the linear dependence (\ref{E:abuvtGamthin}). The first term on the right (\ref{E:oabthin1}) is also $0$ because ${\tilde\Gamma}_u(t_0)$ constant in $u\in [s_0,s_1]$ by assumption (\ref{E:Gammasti}). Hence the path  (\ref{E:uGammathin}) on $\pap$ is $\oab$-horizontal.
 \fbox{QED}

\section{$2$-groups by many names}\label{S:2grp}

In this section we organize known notions, results  and examples in a manner that is useful for our purposes in later sections.  These ideas and related topics on categorical groups  are discussed in the works \cite{Awody, {BaezSchr}, {BaezSchr2}, BarrMack,  Bart, ForBar,{KPT}, Kelly, Macl, {Port },Whit}. 

A category $\mbk$ along with a bifunctor 
\begin{equation}\label{E:Kotimes}\otimes:\mbk\times\mbk\to\mbk \end{equation}
forms a {\em categorical group}   if both $\Obj(\mbk)$ and $\Mor(\mbk)$ are groups under the operaton $\otimes$ (on objects and on morphisms, respectively).
Sometimes it will be more convenient to write $ab$ instead of $a\otimes b$. Being a functor, $\otimes$ carries the identity morphism $(1_a,1_b)$, where $1_x:x\to x$ denotes the identity arrows on $x$,  to the identity morphism $1_{ab}$:
$$1_a\otimes 1_b=1_{a\otimes b}.$$
and so, taking $a$ to be the identity element $e$   in $\Obj(\mbk)$, it follows that $1_e$ is the identity element in the group $\Mor(\mbk)$.

The functoriality of $\otimes$ implies the `exchange law':
\begin{equation}\label{E:exchange}
(f\otimes g)\circ(f'\otimes g')=(f\circ f')\otimes (g\circ g')\end{equation}
for all $f,f',g,g'\in\Mor(\mbk)$ for which the composites $f\circ f'$ and $g\circ g'$ are defined.

The following is a curious but useful consequence of the definition of a categorical group:
 
 \begin{prop}\label{P:2grpHK} Let $\mbk$ be a categorical group, with the group operation written as juxtaposition. 
Then for any morphisms $f:a\to b$ and $h:b\to c$ in $\mbk$, the composition $h\circ f$ can be expressed in terms of the product operation in $\mbk$:
\begin{equation}\label{E:fhcompprod} 
h\circ f =f1_{b^{-1}}h=h1_{b^{-1}}f. \end{equation}
In particular,  
\begin{equation}\label{E:hkkh}
hk=h\circ k=kh\qquad\hbox{if $t(k)=s(h)=e$.} \end{equation}
 \end{prop}

 \noindent\underline{Proof}.  For $f:a\to b$ and $h:b\to c$ we have, on using the exchange law (\ref{E:exchange}),
\begin{equation}\label{E:hicrffhprf}
h\circ f=(1_eh)\circ (f1_{b^{-1}}1_b)=(1_e\circ f1_{b^{-1}})(h\circ 1_b)= f1_{b^{-1}} h.
\end{equation}
Interchanging the order of the multiplication we also have
\begin{equation}\label{E:ficrhfhprf}
h\circ f=(h1_e)\circ (1_b1_{b^{-1}}f)= (h\circ 1_b)(1_e\circ1_{b^{-1}}f)=h1_{b^{-1}}f
\end{equation}
\fbox{QED}
 
 Here is an alternative formulation of the definition of categorical group:

\begin{prop}\label{P:catgrpstix}
If $\mbk$, with operation $\otimes$, is a categorical group then  the source and target maps
$$s, t: \Mor(\mbk)\to \Obj(\mbk)$$
are group homomorphisms, and so is the identity-assigning map
$$ \Obj(\mbk)\to\Mor(\mbk): x\mapsto 1_x.$$
Conversely, if $\mbk$ is a category for which both $\Obj(\mbk)$ and $\Mor(\mbk)$ are groups,  $s$, $t$, and $x\mapsto 1_x$ are homomorphisms, and the exchange law (\ref{E:exchange}) holds then $\mbk$ is a categorical group. 
\end{prop}
\noindent\underline{Proof}. Suppose $\mbk$ is a categorical group. Consider any morphisms $f:a\to b$ and $g:c\to d$ in $\Mor(\mbk)$.  Then, by definition of the product category $\mbk\times \mbk$, we have
the morphism $(f,g)\in\Mor(\mbk\times\mbk)$ running from the domain $(a,c)$ to the codomain $(b,d)$:
$$(f,g):(a,c)\to (b,d)\quad\hbox{in $\Mor(\mbk\times\mbk)$.}$$
 Then, since $\otimes$ is a functor, $f\otimes g$ runs from domain $a\otimes c$ to target $b\otimes d$. Thus,
 $$s(f\otimes g)=s(f)\otimes s(g)\quad\hbox{and}\quad t(f\otimes g)=t(f)\otimes t(g).$$
 Thus $s$ and $t$ are homomorphisms.
 
 Next,   for any objects $x, y\in\Obj(\mbk)$, the identity morphism
 $$(1_{x}, 1_y): (x,y)\to (x,y) $$
 in $ \Mor(\mbk\times\mbk) $
 is mapped by the functor $\otimes$ to the identity morphism
 $$1_{x\otimes y}:x\otimes y\to x\otimes y,$$
and this just means that
 $$1_x\otimes 1_y =1_{x\otimes y}.$$
 Thus $x\mapsto 1_x$ is a group homomorphism.
 
 Conversely, suppose $s$, $t$, $x\mapsto 1_x$ are group homomorphisms and the exchange law  (\ref{E:exchange})  holds (both sides of (\ref{E:exchange})  are meaningful because $s$ and $t$ are homomorphisms). Then (\ref{E:exchange})  says that $\otimes$ maps composites to composites, while $x\mapsto 1_x$  being a homomorphism implies that $1_{x\otimes y}=1_x\otimes 1_y$, and so $\otimes$ is indeed a functor, preserving composition and mapping identities to identities.
 \fbox{QED}
 
 See \cite{ForBar} for a review of $2$-groups, related notions, along with more references (the relation  (\ref{E:fhcompprod}), proved above, is mentioned in \cite{ForBar}.)
 
  The preceding result suggests that when working with groups with more structure, for example with Lie groups, it would be more natural to continue to require that the source, targets, and identity-assigning maps respect the additional structure. Thus  by a {\em categorical Lie group } we mean a category $\mbk$ along with a functor $\otimes$ as above, such that $\Mor(\mbk)$ and $\Obj(\mbk)$ are Lie groups, and the maps $s$, $t$, and $x\mapsto 1_x$ are smooth homomorphisms.

{\bf Example CG1.} For any group $K$ we can construct a categorical group $\mbk_0$ by taking $K$ as the object set and requiring there be a unique morphism $a\to b$ for every $a, b\in K$. At the other extreme we have the discrete categorical group $\mbk_d$ whose object set is $K$ but whose morphisms are only the identity morphisms. $\square$

{\bf Example CG2.} Let  
 \begin{equation}\label{E:hatKK}
 \pi:{\hat K}\to K \end{equation}
 be a surjective homomorphism of groups. 
 We   think of ${\hat K}$ as a `covering group', or a principal bundle  over $K$, with each  fiber $\pi^{-1}(k)$ standing `above' the point $k$. The structure group is 
 $$Z=\ker\pi=\pi^{-1}(e),$$ with $e$ being the identity element of $K$.  For the category $\mbk$ we take the object set to be $K$ itself.  A morphism $a\to b$ is to be thought of as an arrow ${\hat a}\to {\hat b}$, with $\pi({\hat a})=a$ and $\pi({\hat b})=b$ and such that ${\hat a}\to {\hat b}$ is identified with ${\hat a}{\hat k}\to {\hat b}{\hat k}$, for every ${\hat k}\in \ker\pi$.  Thus the discrete categorical group ${\mathbf Z}_d$ acts on the right  on the category $\hat{\mbk}_0$ (notation as in CG1) as follows: an object $z\in Z$ acts on the right on an object ${\hat a}\in {\hat K}$ to produce the product ${\hat a}z$; a morphism $z\to z$ in $\Mor(\mbbz_d)$ acts on a morphism ${\hat a}\to {\hat b}$ in $\Mor(\hat{\mbk}_0)$ to produce the morphism ${\hat a}z\to {\hat b}z$. (See (\ref{E:rightact}) for a precise definition of right action in this context.) We construct a category $\mbk$ which  we view as the quotient ${\hat\mbk}/{\mbbz_d}$. The   object set of $\mbk$ is $K$.  A morphism $f:a\to b$ for $\mbk$ is obtained as follows:  we choose some ${\hat a}\in\pi^{-1}(a)$, ${\hat b}\in\pi^{-1}(b)$ and take $f$ to be the arrow ${\hat a}\to {\hat b}$, identifying this with ${\hat a}z\to {\hat b}z$ for all $z\in Z$. More compactly, $\Mor(\mbk)$ is the quotient $\Mor({\hat\mbk})/\Mor(\mbbz_d)$ where the action of the group $\Mor(\mbbz_d)$ on   $\Mor({\hat\mbk})$ is given by: $({\hat a}\to {\hat b})(z\to z)={\hat a}z\to {\hat b}z$.  In order to make the category $\mbk$ into a  categorical group we define a multiplication on  $\Mor({ \mbk})$ by using the multiplication structure on $\Mor({\hat\mbk})$: for a morphism  $f:a\to b$ given by ${\hat a}\to {\hat b}$ and $g:c\to d$ given by ${\hat c}\to {\hat d}$ we define $fg$ to be
 \begin{equation}\label{E:CG2multip}
 fg=({\hat a}{\hat c}\to {\hat b}{\hat d}).
 \end{equation}
 For this to be well-defined, the right side should not depend on the choices ${\hat x}\in\pi^{-1}(x)$;  if we assume that the subgroup $Z$ is {\em central} in ${\hat K}$,   then:
$${\hat a}z{\hat c}w={\hat a}{\hat c}zw$$
holds for all ${\hat a}, {\hat c}\in {\hat K}$ and $z, w\in Z$, and so the rght side in  (\ref{E:CG2multip}) is determined entirely by $f$ and $g$. $\square$

 A more specific class of examples of categorical    groups is provided by taking $K$ to be any compact Lie group and $\hat K$ its universal covering group.

A group $K$ gives rise to a category $\mbg(K)$ in a natural way: $\mbg(K)$ has just one object $O$,   the morphisms of $\mbg(K)$ are the elements of $K$,  all having  $O$ as both source and target, with composition of morphisms being given by the group operation in $K$.

Next we discuss a concrete form of a categorical group.  For this recall the notion of a {crossed module} $(G, H, \alpha, \tau)$ from (\ref{E:Peiffer}).   
\begin{theorem}\label{T:catgrp2grp} Suppose $\mbg$ is a categorical group, with the group operation written as juxtaposition: $a\otimes b=ab$, and with  $s, t:\Mor(\mbg)\to\Obj(\mbg)$ denoting the source and target maps.  Let $H=\ker s$.  Let
$\tau: H\to G$ and $\alpha:G\to {\rm Aut}(H)$ be the maps specified by
\begin{equation}\label{E:deftaualphacm}
\begin{split}
\tau(\theta) &=t(\theta), \\
\alpha(g)\theta &=1_g\theta 1_{g^{-1}}
 \end{split}\end{equation}
for all $g\in G$ and $\theta\in H$.  
Then  $(G, H, \alpha,\tau)$ is a {crossed module}. 

Conversely, suppose $(G, H, \alpha,\tau)$ is a {crossed module}. Then there is a category $\mbg$   whose object set is $G$ and for which a morphism $h:a\to b$ is an ordered pair
$$(h,a),$$
where $h\in H$ satisfies
 $$\tau(h)a=b,$$
with composition being given by
\begin{equation}\label{E:compinH}
(h_2,b)\circ (h_1,a)=(h_2h_1, a).\end{equation}
Moreover, $\mbg$ is a categorical group, with group operation on $\Obj(\mbg)$ being the one on $G$, and the group operation on $\Mor(\mbg)$ being
\begin{equation}\label{E:multinH}
(h,a)(k,c)=\left(h\alpha(a)(k), ac\right).\end{equation}
\end{theorem}

More fancifully, we have a natural isomorphism between the category of crossed modules and the category of categorical groups, defined in the obvious way.

 Let $(G, H,\alpha,\tau)$ be the crossed module constructed as above from a categorical group $\mbg$, as above.  Let $\mbk$ be the categorical group  constructed from $(G,H,\alpha,\tau)$ according to the prescription in Theorem \ref{T:catgrp2grp}.  The objects of $\mbk$ are just the elements of $G$. A morphism in $\Mor(\mbk)$ is a pair
\begin{equation}\label{E:stphiinv}
\left(s(\phi),\phi 1_{s(\phi)^{-1}}\right), \end{equation}
where $\phi\in\Mor(\mbg)$, so that
\begin{equation}\label{E:stphi1}
\phi 1_{s(\phi)^{-1}}\in H=\ker s.\end{equation}

\noindent\underline{Proof}. Assuming that $\mbg$ is a categorical group we will shows that $(G, H,\alpha,\tau)$ is a crossed module. Let $\tau$ be the target map  $t$ restricted to the subgroup $H$ of $\Mor(\mbg)$:
\begin{equation}\label{E:taugh}\tau:H\to G: h\mapsto \tau(h)= t(h).\end{equation}
Next for $g\in G$ let $\alpha(g):H\to H$ be given  on any $h\in H$ by:
\begin{equation}\label{E:alphahgi}
\alpha(g)h= 1_{g }  h   1_{g^{-1}}.\end{equation}Note that the source of $\alpha(g)h$ is
$$s\bigl(\alpha(g)h\bigr)=g  e  g^{-1}=e,$$
where $e$ is the identity element of $\Obj(\mbg)$. Thus $\alpha(g)h$ is indeed an element of $H$.  Moreover, it is readily checked that $\alpha(g)$ is an automorphism of $H$.

The target of $\alpha(g)h$ is
\begin{equation}\label{E:tauagh}
\tau\bigl(\alpha(g)h\bigr)=t\bigl(\alpha(g)h\bigr)= g   \tau(h)  g^{-1}.\end{equation}

Next, consider $h, h'\in H$, and suppose $\tau(h)=g$ and $\tau(h')=g'$. Then 
\begin{equation}\label{E:ahhh}
\alpha\bigl(\tau(h)\bigr)(h')=1_gh'1_{g^{-1}}=hh'h^{-1},\end{equation}
the last equality of which can be verified by using Proposition \ref{P:2grpHK} (specifically (\ref{E:hkkh})), which implies that
the element $h^{-1}1_g:g\mapsto e$ commutes with $h':e\to g'$.  Equations (\ref{E:tauagh}) and (\ref{E:ahhh}) are exactly the conditions (\ref{E:Peiffer})  that make $(G, H,\alpha,\tau)$ a crossed module.

Before proceeding to the proof of the converse part, we observe that a morphism $f\in\Mor(\mbg)$ is completely specified by its source $a=s(f)$ and by the morphism $$h=f1_{a^{-1}}\in \ker s.$$

Now   suppose $(G,H,\alpha, \tau)$ is a {crossed module}. We construct a category $\mbk$ with object set  $G$.  If $f:a\to b$ is to be a morphism then $  f1_{a^{-1}}$ would be a morphism $e\to ba^{-1} $. Thus, for $a, b\in G$ we take a morphism $a\to b$ to be specified by the source $a$ along with an  element $h\in H$ for which  $\tau(h)=ba^{-1}$.  So we define a morphism $f\in\Mor(\mbk)$ to be an ordered pair $(h,a)\in H\times G$, with source and target given by 
\begin{equation}\label{E:afstbh}
 s(f)=a\qquad\hbox{and}\qquad t(f)=\tau(h)a.
\end{equation}

Thus the source and target maps are
\begin{equation}\label{E:stau}
s(h,a)=a\quad\hbox{and}\quad t(h,a)=\tau(h)a.\end{equation}
To understand what the composition law should be,  recall from (\ref{E:fhcompprod}) that for any morphisms $f:a\to b$ and $g:b\to c$ we have
\begin{equation}
(g\circ f)1_{a^{-1}} = g1_{b^{-1}}f 1_{a^{-1}}.\end{equation}
Thus we can define the composition of morphisms $(h_1,a), (h_2,b)\in H\times G$, with $b=\tau(h_1)a$,  by
\begin{equation}\label{E:compK}
(h_2, b)\circ (h_1, a)=(h_2h_1, a).\end{equation}
  
The identity morphism $1_a$ is then $(e,a)$ because
$$(h_2, a)\circ (e, a)=(h_2, a)\quad\hbox{and}\quad (e, c)\circ (h_1, b)=(h_1, b), $$
where $c=\tau(h_1)b$.
Associativity of composition is clearly valid. Thus $\mbk$ is indeed a category. 
It remains to define a product on $\Mor(\mbk) $ and prove that this product is functorial.

For $f:a\to b$ and $g:c\to d$ we have
$$(fg)1_{(ac)^{-1}}=f1_{a^{-1}} 1_{a }g1_{c^{-1}}1_{a^{-1}}$$
which motivates us to define the product on $\Mor(\mbk)$ by
\begin{equation}\label{E:acbkprod}
(h,a)(k,c)=\left(h \alpha(a)(k), ac\right).\end{equation}
If we identify $H$ and $G$ with  subsets of $\Mor(\mbk)$ through the injections
\begin{equation}\label{E:HGembedMorK}
\begin{split}
H\to\Mor(\mbk) : h &\mapsto (h,e)\\
G\to \Mor(\mbk): a & \mapsto  (e,a)
\end{split}
\end{equation}
 then, using (\ref{E:acbkprod}), the product $ha$ corresponds to the element
 \begin{equation}\label{E:haeah}
ha= (h,e)(e,a)=(h,a).
 \end{equation}
 Hence
  the multiplication operation (\ref{E:acbkprod}) takes the form 
\begin{equation}\label{E:semidirprod}
hakc= h\alpha(a)(k)\,  ac,
\end{equation}
 which means the commutation relation 
 \begin{equation}\label{E:semidcomm}
 ah=\alpha(a)(h)\,a\quad\hbox{for all $a\in G$ and $h\in H$,}
 \end{equation}
 or, equivalently:
  \begin{equation}\label{E:semidcomm2}
 \alpha(a)(h)=aha^{-1}\quad\hbox{for all $a\in G$ and $h\in H$.}
 \end{equation}
 The product (\ref{E:acbkprod}) is a standard semi-direct product of groups, and it is  a straightforward, if lengthy, calculation to verify that $\Mor(\mbk)$ is indeed a group under the operation:
\begin{equation}\label{E:Mormbksemidr}
\Mor(\mbk)=H\rtimes_{\alpha}G.\end{equation}
 It is readily checked that the source and target maps $s$ and $t$ given by (\ref{E:stau}) are homomorphisms, and so is the map
\begin{equation}
\Obj(\mbk)\to \Mor(\mbk):a \mapsto 1_a=(e,a).
 \end{equation}
 Thus $\mbk$ is a categorical group. 

It is straightforward to verify the exchange law  (\ref{E:exchange})  holds. Finally, by construction, $G=\Obj(\mbk)$, and $H$ is  identifiable as the subgroup of $\Mor(\mbk)$  given by $\ker s$.
\fbox{QED}

{\bf Example CG3.}   Let $L_o$ be the set of all backtrack erased piecewise $C^1$ loops  based at a point $o$ in a manifold $B$; under composition, this is clearly a group.    Now the method of example CG2 can be used, with $K=L_o$, to form categorical groups with object group $L_o$. $\square$

 By a {\em $2$-category} $\mbbc_2$ over a category $\mbbc_1$ we mean a   category $\mbbc_2$ whose objects are the arrows of $\mbbc_1$:
 $$\Obj(\mbbc_2)=\Mor(\mbbc_1).$$
  Let $\mbk$ be a categorical group, and $K$ the group formed by $\Obj(\mbk) $. Then $\mbk$ is a $2$-category over $\mbg(K)$, which is the categorical  group whose object group has just one element and whose morphism group is $K$.

   Let $\mbg_1$ and $\mbg_2$ be categorical groups such that $\Obj(\mbg_2)=\Mor(\mbg_1)$.   Let $(G, H,\alpha_1,\tau_1)$ be the    crossed module associated with $\mbg_1$, and $(H\rtimes_{\alpha_1}G, K,\alpha_2,\tau_2)$ the crossed module associated with $\mbg_2$.   We identify $H$ and $G$ naturally with subgroups of $H\rtimes_{\alpha_1}G$, so that each element of this semi-direct product can be expressed as $hg$, with $h\in H$ and $g\in G$. The following computation will be useful later:
   
   \begin{lemma}\label{L:G2H2al2} With notation as above,
   \begin{equation}\label{E:alph2alpha1}
   \alpha_2\Bigl(\alpha_1(g_1^{-1})(h_1^{-1}h)\Bigr)\Bigl(\alpha_2(g_1^{-1})(k)\Bigr)= \alpha_2(g_1^{-1}h_1^{-1}h)(k),
   \end{equation}
   for all $g_1\in G_1$, $h_1, h\in H$, and $k\in K$.
   \end{lemma}
   \noindent\underline{Proof}.  Recall from (\ref{E:semidcomm}) the commutation rule:
\begin{equation}\label{E:gahcommr}
gh=\alpha_1(g)(h)\,g\qquad\hbox{for all $g\in G$ and $h\in H$.}\end{equation}
   Then we have
   \begin{equation}
   \begin{split}
   \alpha_2\Bigl(\alpha_1(g_1^{-1})(h_1^{-1}h)\Bigr)\Bigl(\alpha_2(g_1^{-1})(k)\Bigr) &= \alpha_2\Bigl(\alpha_1(g_1^{-1})(h_1^{-1}h)g_1^{-1}\Bigr)(k)\\
   &\hskip -.5in = \alpha_2(g_1^{-1} h_1^{-1}h)(k)\quad\hbox{(using (\ref{E:gahcommr}) with $g=g_1^{-1}$)},\\
   \end{split}
   \end{equation}
   which proves (\ref{E:alph2alpha1}). \fbox{QED}

\section{Principal categorical bundles}\label{S:categorbundle}

In the traditional topological description, a principal bundle $\pi:P\to B$ is a smooth surjection of manifolds, along with a Lie group $G$ acting freely on $P$ on the right, the action being transitive on each fiber $\pi^{-1}(b)$ (we do not consider local triviality here).  In this section we formulate a categorical description that includes a wider family of geometric objects than is included in the traditional notion of principal bundles. In particular, the framework of categorical bundles  includes bundles over pathspaces with a pair of groups serving as structure groups, one at the level of objects and one at the level of morphisms.  Categorical bundles, in different formalizations, have been studied in the literature (for example in  Baez et al. \cite{BaezDol, BaezLauda, BaezSchr, BaezSchr2, BaezWise}, Bartels \cite{Bart},   Schreiber and Waldorf \cite{SW1,SW2}, Abbaspour and Wagemann \cite{AbbasWag},  and Viennot \cite{Vien}).  Our exploration is distinct and our focus is more on the geometric side than on category theoretic exploration. However, unlike many other works in the area, we do not explore any notions of local triviality nor do we impose any smooth structure on the base categories in the formal definition.

Let $\mbp$ be a category and $\mbbz$ a categorical group.  Denote by $P$ the set of objects of $\mbp$ and by $Z$ the set of objects of $\mbbz$:
$$P=\Obj(\mbp)\quad\hbox{and}\quad Z=\Obj(\mbbz).$$
By a {\em right action} of $\mbbz$ on $\mbp$ we mean a functor
\begin{equation}\label{E:rightact}
\mbp\times\mbbz\to\mbp: (x,g)\mapsto \rho(g)x= \rho(x,g) \end{equation}
that is a right action both at the level of objects and at the level of morphisms; thus, both
\begin{equation}\begin{split}
\Obj(\mbp)\times \Obj(\mbbz) &\to\Obj(\mbp):(A,g)\mapsto \rho(A,g)\\
\Mor(\mbp)\times \Mor(\mbbz) &\to\Mor(\mbp):(F,\phi)\mapsto \rho(F, \phi)
\end{split}\end{equation}
are  right actions. We assume, moreover, that both these actions are {\em free}.

Note that functoriality of $\rho$ implies, in particular, that   
\begin{equation}\label{E:strho}
\begin{split}
s\bigl(\rho(F, \phi)\bigr) &=\rho\bigl(s(F), s(\phi)\bigr)
\\
t\bigl(\rho(F, \phi)\bigr) &=\rho\bigl(t(F), t(\phi)\bigr) 
\end{split}
\end{equation}

By analogy with principal bundles we define a {\em principal categorical bundle} to be a functor $\pi:\mbp\to\mbbb$ along with a right action of a categorical group $\mbbz$ on $\mbp$ satisfying the following conditions:
\begin{itemize}
\item[(i)] $\pi$ is surjective both at the level of objects and at the level of morphisms;
\item[(ii)] the action of $\mbbz$ on $\mbp$ is free on both objects and morphisms;
\item[(iii)] the action of $\Obj(\mbbz)$ on the fiber $\pi^{-1}(b)$ is transitive for each object $b\in\Obj(\mbbb)$, and the action of $\Mor(\mbbz)$ on the fiber $\pi^{-1}(\phi)$ is transitive for each morphism  $\phi\in\Mor(\mbbb)$.
\end{itemize}
 Notice, however, that we are not imposing any form of local triviality.

There is a consistence property of compositions in ${\mathbf P}$ with respect to the action of ${\mathbf Z}$:
 
\begin{lemma}\label{L:comporbits}
 Suppose  $\rho$ is a right action of a categorical group $\mbbz$   on a category $\mbp$, the action being free on both objects and morphisms.  Let $F$ and $H$ be morphisms in $\Mor(\mbp)$, with the composition $H\circ F$ defined, and let $F'$ and $H'$ be morphisms in $\Mor(\mbp)$, with $F'$ being in the $\Mor(\mbbz)$-orbit of $F$ and $H'$ in the $\Mor(\mbbz)$-orbit of $H$, and the composite $H'\circ F'$ also defined; then the composition $H'\circ F'$ is   in the $\Mor(\mbbz)$-orbit of $H\circ F$. 
\end{lemma} 
\noindent\underline{Proof}.  Suppose $F'$ is in the $\Mor(\mbbz)$-orbit of $F$, and $H'$ is in the $\Mor(\mbbz)$-orbit of $H$. Then
\begin{equation}\label{E:F'Fpsi}
F'=F\rho(\phi)\qquad \hbox{and}\qquad H'=H\rho(\psi)\end{equation}
for some $\phi,\psi\in\Mor(\mbbz)$.  
The composability  $H'\circ F'$ implies that the target of $F\rho(\phi)=\rho(F,\phi)$ is the source of $H\rho(\psi)$, and so, by (\ref{E:strho}), 
\begin{equation}\label{E:rhoFHfree}
\rho\bigl(t(F), t(\phi)\bigr)=t\left(\rho(F,\phi)\right)=s\left(\rho(H,\psi)\right)= \rho\bigl(s(H), s(\psi)\bigr). \end{equation}
Now $t(F)=s(H)$ when $H\circ F$ is defined; hence by (\ref{E:rhoFHfree}) and  by freeness of the action $\rho$ on $\Mor(\mbp)$ we conclude that
$$t(\phi)=s(\psi),$$
and so the composite
$\psi\circ\phi$ is defined. Then, using the functoriality of the categorical group action $\rho$, we have
\begin{equation}\label{E:H'F'2}
H'\circ F'=\rho(H,\psi)\circ \rho(F, \phi)=\rho(H\circ F, \psi\circ \phi) 
\end{equation}
which shows that $H'\circ F'$ is in the $\Mor(\mbbz)$-orbit of $H\circ F$. Hence composition is well-defined on the quotient $\Mor(\mbbb)$ as specified in (\ref{E:orbitcomp}). \fbox{QED}

 The consistency property of composition allows us to form a `quotient'  category ${\mathbf P}/{\mathbf Z}$:

\begin{theorem}\label{T:catprinc}  Suppose  $\rho$ is a right action of a categorical group $\mbbz$   on a category $\mbp$, the action being free on both objects and morphisms.   

Let $\mbbb$ be the category whose object set is $B=\Obj(\mbp)/\Obj(\mbbz)$, and whose morphisms are $\Mor(\mbbz)$-orbits of morphisms in $\Mor(\mbp)$, with composition  defined by 
\begin{equation}\label{E:orbitcomp}
[H\circ F]=[H]\circ [F],
\end{equation}
where $[\cdots]$ denotes the $\Mor(\mbbz)$-orbit. Then  
$$\pi:\mbp\to\mbbb$$
taking every object $X\in\Obj(\mbp)$ to the $\Obj(\mbbz)$-orbit of $X$, and every morphism $F$ to its $\Mor(\mbbz)$-orbit is a functor, and, along with the right action of $\mbbz$ on $\mbp$, specifies a principal categorical bundle.
\end{theorem}

\noindent\underline{Proof}. 
Let 
$$B=\Obj(\mbp)/\Obj(\mbbz),$$
and 
$$\pi_P:P\to B$$
the quotient map. 
For $x\in B$ the identity morphism $1_x$ is the orbit of $1_X$, for any $X\in P$ whose orbit is $x$.  Taking $B$ as the set of objects, and morphisms as just described, we have a category $\mbbb$. 

It is clear from the construction of morphisms and compositions in $\mbbb$ that the  quotient association
$$\pi:\mbp\to\mbbb $$
is a functor. Moreover,  again from the construction of  $\mbbb$ in terms of orbits of the action on $\mbbz$, the right action of $\mbbz$ on $\mbp$ is transitive on fibers. By hypothesis, the action is also free. Hence $\pi$, along with the action, specifies a principal categorical bundle.
\fbox{QED}

We can revisit Example CG2 in light of the preceding theorem. For this we take $\mbp$ to be the category ${\hat \mbk}_0$   whose objects are the elements of a given group ${\hat K}$, and whose morphisms are all ordered pairs of elements of ${\hat K}$; let $Z$  be a  subgroup  of ${\hat K}=\Obj({\hat \mbk})$, acting by right multiplication on $\hat K$, and $\mbbz_d$ the categorical group in which the only morphisms are the identities $z\to z $ with $z$ running over $Z$.   Then, as in Example CG2,  there is a category ${\mathbf K}={\hat\mbk}/\mbbz_d$. This provides a categorical principal bundle   ${\hat\mbk}\to \mbk$ with structure categorical group $\mbbz_d$ (for this there is no need to assume $Z$ is central).

{\bf  Example  P1}. A traditional principal $G$-bundle $\pi:P\to B$  generates a categorical principal bundle in the following way. Take $\mbp$ to be the discrete category with object set $P$ (in a discrete category the only morphisms are the identity morphisms), $\mbbb$ to be the discrete category with object set $B$, and $\mbg_d$ to be the categorical group whose object set is $G$ and whose only morphisms are the identity morphisms. Then we have a principal categorical bundle.  

{\bf Example P2.} A more interesting example is obtained again from a principal $G$-bundle $\pi:P\to B$, but with the categorical group (Example CG1) being ${\mathbf G}_0$, whose object set is $G$ and  for which there is a unique morphism $g_1\to g_2$ for every $g_1, g_2\in G$; we denote this morphism by $(g_1,g_2)$.  Take $\mbbb$ to be the category with object set $B$ and with morphisms being backtrack-erased paths on $B$. For $\mbp$ take the category whose object set is $P$ and for which a morphism $p\to q$ is a triple
$$(p,q;\gamma),$$
where $\gamma$ is any backtrack-erased path on $B$ from $\pi(p) $ to $\pi(q)$. Define composition of morphisms in the obvious way, and define a right action of ${\mathbf G}$ on ${\mathbf P}$ by taking it to be the usual right action of $G$ on $P$ at the level of objects and by setting
\begin{equation}\label{E:pqgammag1g2}
(p,q;\gamma)(g_1,g_2)=(pg_1,qg_2;{\gamma}),\end{equation}
for all $p,q\in P$, all backtrack-erased paths $\gamma$ on $B$ from $\pi(p)$ to $\pi(q)$ and all  morphisms $(g_1,g_2)$ in $\mbg_0$.  Defining the projection ${\mathbf P}\to {\mathbf B}$ to be $\pi$ at the level of objects and
$(p,q;\gamma)\mapsto\gamma$ at the level of morphisms yields then a principal categorical bundle.

{\bf Example P3.} Next consider a connection $\ovA$ on a principal $G$-bundle $\pi:P\to B$. Let $\mbbb$ have object set $B$, with arrows being backtrack-erased  paths in $B$. Let $\catmpt^0_{\ovA}$ have objects the points of $P$, with arrows being backtrack erased $\ovA$-horizontal paths. There is then naturally the projection functor $\pi:\mbp\to\mbbb$. The categorical group ${\mathbf G}_d$  whose objects form the group $G$, and whose morphisms are all the identity morphisms, has an obvious right action on $\catmpt^0_{\ovA}$, and thus we have a categorical principal bundle.

In the following section we will explore a more substantive example.

First, however, we explore   the notion of reduction of a bundle. The basic example of a principal bundle is that of the frame bundle ${\rm Fr}_M$ of a manifold $M$: a typical point $p\in {\rm Fr}_M$ is a basis $(u_1,\ldots, u_n)$ of a tangent space $T_mM$ to $M$ at some point $m$, and $\pi:P\to M$ is defined by $\pi(p)=m$. The group $GL(\mbr^n)$ acts on the right on ${\rm Fr}_M$  by transformations of frames:
$$(u_1,\ldots, u_n)\cdot g=[u_1,\ldots, u_n]\left[\begin{matrix}g_{11} & g_{12} &\ldots &g_{1n}\\
g_{21}&g_{22}&\ldots & g_{2n}\\
\vdots &\vdots &\vdots &\ldots &\vdots\\
g_{n1} & g_{n2} &\ldots & g_{nn}\end{matrix}\right].$$
With this structure $\pi:{\rm Fr}_M\to M$ is a principal $GL(\mbr^n)$-bundle.  If $M$ is equipped with a metric, say Riemannian, then we can specialize to orthonormal bases; two such bases are related by an orthogonal matrix. Let be  the subgroup $O(n)$ of orthogonal matrices inside $GL(\mbr^n)$. Then we have the bundle ${\rm OFr}_M\to M$ of orthonormal frames, and this is a principal $O(n)$-bundle. There is the natural `inclusion' map ${\rm OFr}_M\to {\rm Fr}_M$, which preserves the principal bundle structures  in the obvious way.   The general notion here is that of `reduction' of a principal bundle. Let $\pi:P\to M$ be a principal $G$-bundle and $\beta: G_o\to G$ be a homomorphism of Lie groups; then  a {\em reduction} of $\pi$ by $\beta$ to $G_o$ is a principal $G_0$-bundle $\pi_o:P_o\to M$ along with a smooth map
$$f:P_o\to P$$
that maps each fiber $\pi_o^{-1}(m)$ into the fiber $\pi^{-1}(m)$ and $f(pg)=f(p)\beta(g)$, for all $p\in P_o$ and $g\in G_o$.

We can define an analogous notion for principal  categorical  bundles. Suppose $\pi:{\mathbf P}\to {\mathbf B}$ is a principal categorical bundle with structure group ${\mathbb G}$, a categorical group. Next suppose ${\mathbf G}_o$ is a categorical group and
$$\beta: {\mathbf G}_o\to {\mathbf G}$$
is a functor that is a homomorphism on objects as well as morphisms. Then by a {\em reduction} of $\pi:{\mathbf P}\to {\mathbf B}$ by $\beta$ we mean a principal categorical bundle $\pi_o:{\mathbf P}_o\to {\mathbf B}$, with structure group ${\mathbf G}_o$, along with a functor
$$f: {\mathbf P}_o\to {\mathbf P}$$
that is fiber preserving both on objects and morphisms and satisfies
$$f(pg)=f(p)\beta(g)$$
both when $(p,g)\in \Obj({\mathbf P}_0)\times\Obj({\mathbf G}_0)$ and when $(p,g)\in \Mor({\mathbf P}_0)\times\Mor({\mathbf G}_0)$. We will see an example of such a reduction in Proposition \ref{P:reducex}  in the next section.

\section{A decorated principal categorical bundle}\label{S:pcatpath}

In this section we introduce a new notion, that of a decorating a principal categorical bundle. This process   allows one to use an ordinary principal $G$-bundle $\pi:P\to M$ with connection and a crossed module $(G, H, \alpha, \tau)$  to obtain a principal categorical bundle, over a space of paths on $M$, with structure  categorical group corresponding to $(G, H, \alpha, \tau)$.  We will continue to explore this construction in later sections, in the context of pathspaces and in the abstract.  

Let $\mbk$ be a categorical group and $(G, H, \alpha,\tau)$ the corresponding crossed module. Recall that each morphism $\phi\in\Mor(\mbk)$ can be identified with the ordered pair
$(h,g)\in H\rtimes_{\alpha}G$ where $g=s(\phi)$ and $h={\phi}1_{s(\phi)^{-1}}\in H$, the subgroup of $\Mor(\mbk)$ consisting of morphisms that have source the identity element $e$ in $G=\Obj(\mbk)$. The composition law for morphisms translates to the product in $H$, as explained in  (\ref{E:compK}).

  Let $\catmbt$ be the category whose objects are 
the points of $M$ and whose morphisms are the backtrack erased piecewise $C^1$ paths on $M$. (Note that by `backtrack-erased path' we mean an equivalence class of paths that are backtrack   equivalent to each other.) Instead of working with piecewise $C^1$ paths we could also work with $C^\infty$ paths that are constant near the initial and final points; the latter condition ensures that composites of such paths are $C^\infty$.

We construct a category  $\catmpt_{\ovA}^{\rm dec}$  by decorating the morphisms of  the category in Example P3 in section \ref{S:categorbundle} with elements of the group $H$. Specifically, we define the category  $\catmpt_{\ovA}$ as follows:  (i)  the object set is  $P$; (ii)  a morphism $f$ is a pair  $( {\tilde\gamma},h)$, where ${\tilde\gamma}$ is a  piecewise $C^1$ backtrack-erased $\ovA$-horizontal path on $P$,   and $h\in H$ (corresponding to a morphism $\phi\in\Mor(\mbk)$ with $s(\phi)=e$), source and targets being specified by
\begin{equation}\label{E:stgphi}
s( {\tilde\gamma},h)=s({\tilde\gamma})\quad\hbox{and}\quad t( {\tilde\gamma},h)= t({\tilde\gamma})\tau(h^{-1}).
\end{equation}
Define composition of morphisms by
\begin{equation}\label{E:compmor}
({\tilde\gamma}_2,h_2 )\circ ( {\tilde\gamma}_1,h_1)=({\tilde\gamma}_3, h_2h_1),
\end{equation}
where ${\tilde\gamma}_3$ is the composite  of the path ${\tilde\gamma}_1$ with the right translate ${\tilde\gamma_2}\tau(h_1)$ (so that the final and initial points match correctly):
\begin{equation}\label{E:tilgamcomp}
{\tilde\gamma}_3={\tilde\gamma_2}\tau(h_1)\circ {\tilde\gamma}_1.\end{equation}
Note that
\begin{equation}\label{E:targetcheck}
t({\tilde\gamma}_3,h_2h_1)=t({\tilde\gamma}_2)\tau(h_1) \tau(h_2h_1)^{-1}=t({\tilde\gamma}_2,h_2),
\end{equation}
and the corresponding result for the sources is clear.

Then $ \catmpt_{\ovA}^{\rm dec}$ is a category, and there is the functor
\begin{equation}\label{E:pidec}
\pi:\catmpt_{\ovA}^{\rm dec}\to \catmbt \end{equation}
 that associates to each object $p\in P$ the object $\pi(p)\in M$, and to each morphism $({\tilde\gamma}, h)$  associates the backtrack erased path $\pi\circ{\tilde\gamma}$.

Now we define a right action  $\rho$ of the categorical group $\mbk$ on $\catmpt_{\ovA}^{\rm dec}$. At the level of objects this is simply the usual right action of $G$ on $P$. For morphisms we define the action as follows. Recall that a morphism in $\mbk$ has the form $(h_1, g_1)$, where 
$$s(h_1, g_1)=g_1,\quad t(h_1, g_1)= \tau(h_1)g_1.$$

Then for $({\tilde\gamma},h),\in \Mor(\catmpt_{\ovA})$ and $ (h_1,g_1)\in\Mor(\mbk)$, we define
\begin{equation}\label{E:morrightactpath}
\rho\left(({\tilde\gamma},h),  (h_1,g_1)\right)= ({\tilde\gamma},h)\cdot (g_1,h_1)
=\left( {\tilde\gamma}g_1,   \alpha(g_1^{-1})(h_1^{-1}h)\right),
\end{equation}
where, on the right, ${\tilde\gamma}g_1$ is the usual right-translate of ${\tilde\gamma}$ by $g_1$.

 The following computation shows that (\ref{E:morrightactpath}) does define a right action:
\begin{equation}\label{E:rightactverfy1}
\begin{split}
\left[( {\tilde\gamma},h)\cdot (h_1,g_1)\right]\cdot (h_2,g_2) &=\left({\tilde\gamma}g_1, \alpha(g_1^{-1})(h_1^{-1}h) \right)\cdot (h_2,g_2) \\
&=\left({\tilde\gamma}g_1g_2, \,\alpha(g_2^{-1})\left( h_2^{-1}\alpha(g_1^{-1}) (h_1^{-1}h)\right) \right)\\
&=\left({\tilde\gamma}g_1g_2, \alpha((g_1g_2)^{-1})(\alpha(g_1)(h_2^{-1})h_1^{-1}h)]\right)\\
&= (  {\tilde\gamma}, h)\cdot \left(  h_1 \alpha(g_1)(h_2), g_1g_2  \right)\\
&= ({\tilde\gamma}, h)\cdot\left[(h_1,g_1 )\cdot (h_2,g_2)\right] \qquad \hbox{(using (\ref{E:acbkprod})).}
\end{split}\end{equation}

From the definition   (\ref{E:morrightactpath}) we see directly that    $({\tilde\gamma}, h)\cdot (h_1, g_1)$ is equal to $ ({\tilde\gamma}, h)$ if and only if $g_1=e$ and $h_1=e$, and so  the action is free. It is also readily checked that the action is transitive on fibers: suppose $( {\tilde\gamma}_0, h_0)$ and $( {\tilde\gamma}, h)$ project to $\gamma$; then   ${\tilde\gamma}_0={\tilde\gamma}g_1$ for a unique $g_1\in G$, and then on taking 
$$h_1= h \alpha(g_1)(h_0^{-1}),$$
we have
$$({\tilde\gamma}_0, h_0)=({\tilde\gamma}, h)\cdot (h_1, g_1).$$

For computational purposes we switch back to the regular morphism notation:  every morphism in $\catmpt_{\ovA}^{\rm dec}$ is of the form 
$$({\tilde\gamma},  \theta)\in  \Mor(\catmpt_{\ovA})\times \ker s  \subset \Mor(\catmpt_{\ovA})\times  \Mor(\mbk)$$
where $s(\theta)=e$,  the identity in $\Obj(\mbk)$.  The source and targets are
\begin{equation}\label{E:sourcetarggamphith}
s({\tilde\gamma},  \theta)=s({\tilde\gamma})\quad\hbox{and}\quad t({\tilde\gamma}, \theta)=t({\tilde\gamma})t(\theta)^{-1}.
\end{equation}

We will now check that   the  expression for the action (\ref{E:morrightactpath}) in the morphism notation is
\begin{equation}\label{E:morrightactpath2}
({\tilde\gamma}, \theta)\cdot\phi= ({\tilde\gamma}s(\phi), \phi ^{-1}\theta 1_{s(\phi)} ) 
\end{equation} 
where $ \phi ^{-1}:a^{-1}\to b^{-1}$ is the multiplicative inverse of $\phi:a\to b\in \Mor(\mbk)$ (we rarely need to use the compositional inverse $b\to a$  and avoid introducing a new notation relying on the context to make the intended meaning clear). To make the comparison we take 
$$\theta= h,$$
 and we recall from (\ref{E:stphiinv}) that a morphism $\phi$ of $\Mor(\mbk) $ is of the form
\begin{equation}\label{E:mrgstphi}
( h_1,g_1)=\Bigl(\phi 1_{s(\phi)^{-1}}, s(\phi)\Bigr).
\end{equation}
 It is clear that the  first component on the right side in (\ref{E:morrightactpath2})   corresponds to the first component on the right side in (\ref{E:morrightactpath}). For the second component we have
\begin{equation}\label{E:comparison}
\begin{split}
\alpha(g_1)^{-1}(h_1^{-1}h) &=\alpha(s(\phi)^{-1})\Bigl(1_{s(\phi) }\phi^{-1} \theta  \Bigr)\\
&=1_{s(\phi)^{-1}}\Bigl(1_{s(\phi) }\phi^{-1} \theta  \Bigr)1_{s(\phi) }\quad\hbox{(using (\ref{E:deftaualphacm})})\\
&=  \phi^{-1} \theta  1_{s(\phi) },
\end{split}
\end{equation}
which is exactly the first component on the right in (\ref{E:morrightactpath2}). Thus  the right action given by (\ref{E:morrightactpath2}) is the same as that given by (\ref{E:morrightactpath}).

The source and target assignments behave functorially:
\begin{equation}\label{E:stfunctphithet}
\begin{split}
s\left(({\tilde\gamma},  \theta )\cdot\phi\right) &= s\left({\tilde\gamma}s(\phi),  \phi^{-1}\theta 1_{s(\phi)},  \right)=s({\tilde\gamma})s(\phi)=s({\tilde\gamma}, \theta)s(\phi)
\\
 t\left(({\tilde\gamma}, \theta)\cdot\phi\right)&=t\left({\tilde\gamma}s(\phi), \phi^{-1}\theta 1_{s(\phi)}\right)\\
 &=t({\tilde\gamma})s(\phi)t\bigl(\phi^{-1}\theta 1_{s(\phi)}\bigr)^{-1} =t({\tilde\gamma})t(\theta)^{-1}t(\phi) \\
 &=t({\tilde\gamma}, \theta)t(\phi). \end{split}
\end{equation}

Next we check that (\ref{E:morrightactpath2}) does specify a right action:
\begin{equation}\begin{split}
\bigl(({\tilde\gamma}, \theta)\cdot \phi_1\bigr)\cdot\phi_2&= ({\tilde\gamma}s(\phi_1), \phi_1^{-1}\theta1_{s(\phi_1)}\bigr)\cdot \phi_2\\
&=\left({\tilde\gamma}s(\phi_1)s(\phi_2),  \phi_2^{-1}\phi_1^{-1}\theta  1_{s(\phi_1) }1_{s(\phi_2) }  \right)\\
&=({\tilde\gamma}, \theta)\cdot (\phi_1\phi_2).
\end{split}\end{equation}
(We have already proved this in (\ref{E:rightactverfy1}) in terms of the crossed module $(G,H,\alpha,t)$.)

The composition law (\ref{E:compmor}) reads
\begin{equation}\label{E:complawpathdec}
({\tilde\gamma}_2, \theta_2)\circ ({\tilde\gamma}_1, \theta_1)=
\left( {\tilde\gamma}_2t(\theta_1) \circ {\tilde\gamma}_1, \theta_2\theta_1\right),
\end{equation}
whose target is clearly the same as the target of $({\tilde\gamma}_2,\theta_2)$ (note that the composite  on the right in (\ref{E:complawpathdec}) is defined). We check associativity:
\begin{equation}\label{E:assoccompos}\begin{split}
({\tilde\gamma}_3, \theta_3)\circ\left( ({\tilde\gamma}_2, \theta_2)\circ ({\tilde\gamma}_1, \theta_1)\right) &= ({\tilde\gamma}_3, \theta_3)\circ \left({\tilde\gamma}_2t(\theta_1) \circ {\tilde\gamma}_1, \theta_2\theta_1\right)\\
&=\left({\tilde\gamma}_3t(\theta_2\theta_1)\circ {\tilde\gamma}_2t(\theta_1)  \circ {\tilde\gamma}_1\, ,\, \theta_3\theta_2\theta_1  \right),
\end{split}\end{equation}
while
\begin{equation}\begin{split}
\left(({\tilde\gamma}_3, \theta_3)\circ ( {\tilde\gamma}_2, \theta_2)\right)\circ ({\tilde\gamma}_1, \theta_1) &= ({\tilde\gamma}_3 t(\theta_2) \circ{\tilde\gamma}_2, \theta_3\theta_2) \circ ({\tilde\gamma}_1, \theta_1)\\
&=\left(  \bigl({\tilde\gamma}_3 t(\theta_2) \circ{\tilde\gamma}_2\bigr)t(\theta_1) \circ {\tilde\gamma}_1 \, ,\, \theta_3\theta_2\theta_1\right),
\end{split}\end{equation}
in agreement wth the last expression in (\ref{E:assoccompos}). The existence of identity morphisms is readily verified.

We can now check functoriality of the right action by examining 
$$\left(({\tilde\gamma}
_2,  \theta_2)\circ ({\tilde\gamma}_1,  \theta_1)\right)\cdot  (\phi_2\circ\phi_1),$$
wherein note that the composablity $\phi_2\circ\phi_1$ means that
\begin{equation}\label{E:stphi12}
s(\phi_2)=t(\phi_1).\end{equation} 

 Composing  and then acting produces:
\begin{equation}\label{E:morfunctcat}
\begin{split}
\left(({\tilde\gamma}
_2,  \theta_2)\circ ({\tilde\gamma}_1,  \theta_1)\right)\cdot  (\phi_2\circ\phi_1) &=
\left( {\tilde\gamma}_2t(\theta_1)\circ {\tilde\gamma}_1,  \theta_2\theta_1\right)\cdot (\phi_2\circ\phi_1) \\
&=\left(\bigl( {\tilde\gamma}_2t(\theta_1 )\circ {\tilde\gamma}_1 \bigr)s(\phi_ 1), (\phi_2\circ\phi_1)^{-1} \theta_2\theta_1 1_{s(\phi_1)}\right),
\end{split}\end{equation}
while first acting and then composing produces:
\begin{equation}\label{E:morfunctcat2}
\begin{split}
({\tilde\gamma}
_2,  \theta_2)\cdot\phi_2\,\,\circ\,\, ({\tilde\gamma}
_1,  \theta_1)\cdot\phi_1 &=\\
\left({\tilde\gamma}_2s(\phi_2), \phi_2^{-1}\theta_21_{s(\phi_2)}\right)\circ \left({\tilde\gamma}_1s(\phi_1), \phi_1^{-1}\theta_11_{s(\phi_1)} \right) &\\
&\hskip -2.5in = \left({\tilde\gamma}_2\underbrace{s(\phi_2)[t(\phi_1)^{-1}t(\theta_1)s(\phi_1)] }_{t(\theta_1)s(\phi_1)}\circ {\tilde\gamma}_1s(\phi_1)\,,\, \phi_2^{-1}\theta_21_{s(\phi_2)}\phi_1^{-1}\theta_11_{s(\phi_1)} \right)
\end{split}\end{equation}
(where in the last step we used $s(\phi_2)=t(\phi_1)$) 
which clearly agrees, in the second entry, with the right side of (\ref{E:morfunctcat}); agreement in the first entry follows using the identities:
\begin{equation}\label{E:morfunctcat3}\theta_2^{-1}(\phi_2\circ\phi_1)\, \stackrel{\rm (\ref{E:fhcompprod} )}{=}\,  \theta_2^{-1} \phi_11_{t(\phi_1)^{-1}}\,\phi_2\stackrel{\rm (\ref{E:hkkh})}{=}\phi_11_{t(\phi_1)^{-1}}\theta_2^{-1}\phi_2.
\end{equation}

We have thus proved:

\begin{theorem}\label{T:pathbundlecat}  Let $(G, H,\alpha,\tau)$ be a  Lie {crossed module} corresponding to a categorical group $\mbk$, $\ovA$ a connection form on a principal $G$-bundle $\pi:P\to M$. Let  $\catmbt$ and $\catmpt_{\ovA}^{\rm dec}$ be the categories constructed above, and $\pi:\catmpt_{\ovA}^{\rm dec}\to \catmbt$ the functor given in (\ref{E:pidec}). Then $\pi:\catmpt_{\ovA}^{\rm dec}\to \catmbt $, along with the right action of $\mbk$ on $\mbp$ defined in (\ref{E:morrightactpath}) is a principal categorical bundle.  
\end{theorem}

In section  \ref{S:catcon} (equation (\ref{E:defmbpabsdecbun})) we will construct a categorical principal bundle
$$\pi:\mbp_{\ovAb}^{\rm dec}\to\mbbb$$
starting from a principal bundle $\mbp \to\mbbb$ and some additional data. This will generalize the specific construction of Theorem \ref{T:pathbundlecat} to provide a `decorated' version of a given categorical principal bundle $\mbp \to\mbbb$.

Now recall that from the original principal bundle $\pi:P\to M$ and the connection $\ovA$ we also have an undecorated principal categorical bundle $\catmpt^0_{\ovA}$ for which the categorical structure group  ${\mathbf G}_d$ has object set $G$ and all the morphisms are the identity morphisms; the morphisms of $\catmpt^0_{\ovA}$ are simply the $\ovA$-horizontal backtrack-erased paths in $P$.   The following result, whose proof is quite clear, is worth noting:
 
\begin{prop}\label{P:reducex} Let  ${\mathbf G}$ be  a categorical group corresponding   to  a Lie crossed module  $(G, H, \alpha, \tau)$, and ${\mathbf G}_d$ be the categorical group whose object set is $G$ and all of whose morphisms are the identity morphisms.   Let
$$R_d:{\mathbf G}_d\to {\mathbf G}$$
be the identity map on objects and the  inclusion map on morphisms. Let $\ovA$ be a connection on a principal $G$-bundle $\pi:P\to M$, and let  $\catmpt_{\ovA}^{\rm dec}$ and $\catmpt^0_{\ovA}$ be  the principal categorical bundles  described above. Consider the association  
\begin{equation}\label{E:redfunct}
R: \catmpt^0_{\ovA}\to \catmpt_{\ovA}^{\rm dec}
\end{equation}
that, at the level of objects, is the identity map on $P$ and for morphisms is given by
$$R({\tilde\gamma})=(e_H, {\tilde\gamma}),$$
where $e_H$ is the identity element in $H$.  Then $R$ is a reduction by $R_d$, in the sense that $R$ maps each fiber into itself, both on objects and on morphisms, and 
$$R({\tilde\gamma}\phi)=R({\tilde\gamma})R_d(\phi)$$
for  all ${\tilde\gamma}\in \Mor\left(\catmpt^0_{\ovA}\right)$ and $\phi\in \Mor({\mathbf G}_d)$.
\end{prop}

\section{Categorical connections}\label{S:catcon}

By a {\em connection} $\mbba$ on a principal categorical  bundle $\pi:\mbp\to\mbbb$, with structure categorical group $\mbk$, we mean a prescription for lifting morphisms in $\mbbb$ to morphisms in $\mbp$.  More specifically, for each $p\in \Obj(\mbk)$ and morphism $ \gamma\in \Mor(\mbbb)$, with source $\pi(p)$,  a connection specifies a morphism ${ \gamma}^{\rm hor}_p:p\to q$, for some $q$ with $\pi(q)=t(\gamma)$; we call the lift ${ \gamma}^{\rm hor}_p$   the {\em lift} of $\gamma$ {\em through $p$}.  Of course,  we require
$$\pi\bigl({\gamma}^{\rm hor}_p\bigr)=\gamma,$$
and the lifting   should be functorial: if $\zeta$ is also a morphism with source $\pi(q)$ then the lift of ${\zeta}\circ{\gamma}$ through $p$ is
\begin{equation}\label{E:horliftfunct}
{ \zeta}^{\rm hor}_q\circ { \gamma}^{\rm hor}_p.\end{equation}
Furthermore, we require that  a `rigid vertical motion' of a horizontal morphism  produces a horizontal morphism:  
\begin{equation}\label{E:horliftrightact}
\gamma^{\rm hor}_{p}1_g
=\gamma^{\rm hor}_{pg} \end{equation}
for all $g\in \Obj(\mbg)$. Identifying morphisms of $\mbk$ with pairs  $(h,g_1)\in H\rtimes_{\alpha}G$, this condition reads
\begin{equation}\label{E:horliftrightact2}
\gamma^{\rm hor}_{pg}=\gamma^{\rm hor}_{p}\cdot (e,g)
\end{equation}
The morphism $\gamma^{\rm hor}_{p}$ will be called the {\em $\mbba$-horizontal lift} of $\gamma$ through $p$. Such morphisms will be called {\em $\mbba$-horizontal morphisms.}

{\bf Example CC1.} The simplest example is the usual connection on a principal bundle. For this, let $\ovA$ be a connection on a principal bundle $\pi:P\to B$. Let $\mbp_0$ be the category whose object set is $P$ and for which a morphism $p\to q$ is a triple $(p,q;\gamma)$, where $\gamma$ is any backtrack-erased path from $\pi(p)$ to $\pi(q)$ on $B$. We take $\mbbb$ to be the category with object set $B$ and morphisms being backtrack-erased paths. Let ${\mathbf G}_0$ have object set $G$, and have a unique morphism $g_1\to g_2$, denoted $(g_1,g_2)$, for every $g_1, g_2\in G$. The action of ${\mathbf G}_0$ on ${\mathbf P}_0$ is as described in (\ref{E:pqgammag1g2}), and then we have a principal categorical bundle $\pi:{\mathbf P}_0\to {\mathbf B}$. For any $p\in P$ and any parametrized piecewise $C^1$ path $\gamma$ on $B$, corresponding to a morphism of $\mbbb$ with source $\pi(p)$, let ${\tilde\gamma}_p$ be the  morphism of $\mbp_0$ specified by the path $\gamma$, the source point $p$ and the target is the point of $P$ obtained by $\ovA$-parallel-transporting $p$ along $\gamma$ to its end. Note that this target, and hence ${\tilde\gamma}_p$,  is determined by    the backtrack-erased form of $\gamma$ along with $p$. Thus, we have a connection $\overline{{\mathbf A}_0}$ on the principal categorical bundle ${\mathbf P}_0\to {\mathbf B}$.

{\bf Example CC2.} Let $\ovA$ be a connection form on a principal $G$-bundle $\pi:P\to M$, and ${\mathbf K}$ a categorical Lie group with {Lie crossed module} $(G, H,\alpha,\tau)$. We now describe a categorical connection on the principal categorical bundle $\pi:\catmpt_{\ovA}^{\rm dec}\to \catmbt $. Let $p\in P$ and $\gamma$ a parametrized piecewise $C^1$ path on $M$ with initial point $\pi(p)$. Define the lift   of the morphism $\gamma$, with source $p$,  to be 
$$(e_H, {\tilde\gamma}_p),$$
where ${\tilde\gamma}_p$ is the backtrack-erased form of the $\ovA$-horizontal lift of $\gamma$ initiating at $p$, and $e_H$ is the identity element in $H$. 

We have the following more general construction of a connection on $\catmpt_{\ovA}^{\rm dec}$. Suppose $(G, H, \alpha,\tau)$ is a Lie crossed module corresponding to a categorical Lie group $\mbk$,   let $\ovA$ be a connection on a principal $G$-bundle $P\to M$, and let  $C$ be an $L(H)$-valued  $C^\infty$ $1$-form on $P$ that is $\alpha$-equivariant in the sense that 
$$C_{pg}(vg)=\alpha(g^{-1})C_p(v)$$
 for all $p\in P$, $g\in G$, and $v\in T_pP$, where $pg=R_g(p)$ is from the right action of $G$ on $P$, and $vg=R_g'(p)v$.

\begin{theorem}\label{T:connectpapbt}  We use notation and hypotheses as above. For any backtrack-erased path $\gamma:[t_0,t_1]\to M$, and   any point $u$ on the fiber over   $s(\gamma)=\gamma(t_0)$, define
\begin{equation}\label{E:hhrliftB0}
\gamma^{\rm hor}_u=
\left({\tilde\gamma}_u, h_u(\gamma)\right),\end{equation}
where ${\tilde\gamma}_u$ is the $\ovA$-horizontal lift of $\gamma$ through $u$ and $h_u(\gamma)$ is the final point of the path $[t_0,t_1]\to H:t\mapsto h(t) $  satisfying the differential equation
\begin{equation}\label{E:htB0}
 h(t)^{-1}h'(t)=-C\bigl({\tilde\gamma}_u'(t)\bigr) \end{equation}
 at all $t\in [t_0,t_1]$ where $\gamma'(t)$ exists,
with initial point $h(t_0)=e$, the identity in $H$. Then this lifting process defines a connection on the categorical principal bundle $\pi:\catmpt_{\ovA}^{\rm dec}\to \catmbt $.\end{theorem}

\noindent\underline{Proof}.  Replacing $u$ by $ug$ in (\ref{E:htB0}), for any $g\in G$, is equal to applying $\alpha(g^{-1})$, and so, by uniqueness of solution of the differential equation, the solution $h(t)$ gets replaced by $\alpha(g^{-1})h(t)$ (recall that $\alpha(g)$ is an automorphism of $H$ and $(g,h)\mapsto \alpha(g^{-1},h)$ is smooth). Hence
\begin{equation}\label{E:hug1}
h_{ug}(\gamma)=\alpha(g^{-1})h_u(\gamma),\end{equation}
and so
\begin{equation}\label{E:hugamtild}
\left({\tilde\gamma}_u, h_{u}(\gamma)\right)\cdot(e,g)=\left({\tilde\gamma}_ug, \alpha(g^{-1})   h_{u}(\gamma)\right)=\left({\tilde\gamma}_{ug}, h_{ug}(\gamma)\right),
\end{equation}
confirming the property (\ref{E:horliftrightact2}).

Now consider backtrack-erased paths $\gamma_1$ and $\gamma_2$ on $M$ for which the composite $\gamma_2\circ\gamma_1$ exists, and let $u$ be in the fiber above $s(\gamma_1)$. Then
$$\gamma^{\rm hor}_u=\left({\tilde\gamma}_{1,u}, h_1\right),$$
where
\begin{equation}\label{E:defh1hugam1}
h_1=h_u(\gamma_1),\end{equation}
has target
\begin{equation}\label{E:vprinetaarg}
v'=t\left(\gamma^{\rm hor}_{1,u}\right)= v \tau(h_1)^{-1},\end{equation}
on using (\ref{E:stgphi}), where
\begin{equation}\label{E:defvendpt}
v\stackrel{\rm def}{=}  t\left({\tilde\gamma}_{1,u}\right)
\end{equation}
is the endpoint of ${\tilde\gamma}_{1,u}$.  The $\ovA$-horizontal lift of $\gamma_2$ with initial point $v'$ is
$${\tilde\gamma}_{2,v'}={\tilde\gamma}_{2,v}\tau(h_1)^{-1}.$$
Next note that
\begin{equation}
\gamma^{\rm hor}_{2,v}=\left({\tilde\gamma}_{2,v}, h_2 \right)\end{equation}
where
\begin{equation}\label{E:defh1hvgam1}
h_2=h_v(\gamma_2).\end{equation}
Then using (\ref{E:vprinetaarg}) , (\ref{E:hugamtild}) and the formula for the right action given in (\ref{E:morrightactpath}), we have
\begin{equation}\begin{split}
\gamma^{\rm hor}_{2,v'} &=\gamma^{\rm hor}_{2,v}\cdot (e, \tau(h_1)^{-1}) \\
&=\left({\tilde\gamma}_{2,v}\tau(h_1)^{-1}\,,\, {h'_2}\right),
\end{split}\end{equation}
where
\begin{equation}\label{E:defh'2}
h'_2=\alpha(\tau(h_1))  h_2.
\end{equation}
Composing with the lift of $\gamma_1$ we have
\begin{equation}\label{Ehorh1h2}
\gamma^{\rm hor}_{2,v'}\circ \gamma^{\rm hor}_{1,u}=\left({\tilde\gamma}_{2,v}\circ {\tilde\gamma}_{1,u},  h'_2h_1\right).
\end{equation}
The second term here is clearly
$${\tilde\gamma}_{2,v}\circ {\tilde\gamma}_{1,u}={\widetilde {(\gamma_2\circ\gamma_1)}}_u.$$
The first term is
$$h'_2h_1= h_1h_2h_1^{-1}\,h_1=h_1h_2.$$
From the differential equation (\ref{E:htB0}) we know that any solution when left multiplied by a constant term in $H$ is again a solution, just with a new initial condition. 
Let $\gamma_1$ have parameter domain $[t_0,t_1]$ and $\gamma_2$ have $[t_1,t_2]$. Hence $h_1h_2$ is  the terminal value $a(t_2)$ of the solution $t\mapsto a(t)\in H$ of 
$$a(t)^{-1}a'(t)=-C\bigl({\tilde\gamma}_{2,v}'(t)\bigr) \quad\hbox{for $t\in [t_1,t_2]$,}  $$
with initial value $h_1$. But $h_1$ itself is the terminal value of the solution of 
$$a(t)^{-1}a'(t)=-C\bigl({\tilde\gamma}_{1,u}'(t)\bigr)  \quad\hbox{for $t\in [t_0,t_1]$,} $$
with initial value $a(t_0)=e\in H$. Thus, fitting these two differential equations into one, we see that $h_1h_2$ is the terminal value $a(t_2)$ of the solution $a(\cdot)$ of the equation
$$a(t)^{-1}a'(t)=-C\bigl({\tilde\gamma}_{u}'(t)\bigr) \quad\hbox{for $t\in [t_0,t_2]$,}  $$
where $\gamma=\gamma_2\circ\gamma_1$. Here $\gamma$, its $\ovA$-horizontal lift $\tilde\gamma$, and $a(\cdot)$ are piecewise $C^1$. 

Hence
\begin{equation}
h'_2h_1= h_1h_2=a(t_2)=h_u(\gamma_2\circ\gamma_1).
\end{equation}
Recalling the definitions of $h_1$ and $h_2$ from (\ref{E:defh1hugam1}) and (\ref{E:defh1hvgam1}),  let us write this explictly (for future reference and use) as
\begin{equation}\label{E:huhom}
 h_u(\gamma_2\circ\gamma_1)= h_u(\gamma_1)h_v(\gamma_2).
\end{equation}
(We will prove a  general version of this observation in Proposition \ref{P:ptmulti}.)
Returning to (\ref{Ehorh1h2}) we conclude that
\begin{equation}\label{Ehorh1h3}
\gamma^{\rm hor}_{2,v'}\circ \gamma^{\rm hor}_{1,u}=\left({\tilde\gamma}_{u},  h_u(\gamma)\right)=\gamma^{\rm hor}_u,
\end{equation}
where again $\gamma=\gamma_2\circ\gamma_1$.  

Thus
$\gamma\mapsto \gamma^{\rm hor}_u$  takes composites to composites. \fbox{QED}

A special case of the preceding construction is obtained by taking
$C$ of the type
$$ (d\Phi)\Phi^{-1},$$
for some smooth $\alpha$-equivariant function $\Phi:P\to H$. In this case the horizontal lift is given by
\begin{equation}\label{E:horliftPhi}
h^{\rm hor}_u(\gamma)= \left({\tilde\gamma}_u, \Phi(v)\Phi(u)^{-1}\right),\end{equation}
with usual notation, writing $v$ for the endpoint of ${\tilde\gamma}_u$.

Let $ {\mathbf A}_0$ be a connection on   a principal  categorical bundle $\pi:\mbp \to\mbbb $ with structure categorical group $\mbg_0 $.  We will now describe an abstract form of the construction used for $\catmpt_{\ovA}^{\rm dec}$ in Theorem \ref{T:connectpapbt}. 

Consider   the category $\mbp_{{\mathbf A}_0}$ whose   object set is $\Obj(\mbp)$ and morphisms are ${\mathbf A}_0$-horizontal morphisms.  The discrete categorical group ${\mbg_0}_d$, whose object group  is $G=\Obj(\mbg_0)$, acts on the right on $\mbp_{{\mathbf A}_0}$ by restriction of the original action of $\mbg_0$ on $\mbp$ as stated in (\ref{E:horliftrightact}). Thus, restricting the projection functor $\pi$ in the obvious way, we see that
\begin{equation}\label{E:piPAbar}
\pi:\mbp_{{\mathbf A}_0}\to\mbbb
\end{equation}
is a principal categorical bundle with structure group ${\mbg_0}_d$. Now consider a categorical Lie group $\mbg_1$ whose object group is $G$. We will construct a principal categorical bundle with structure categorical group $\mbg_1$ by suitably decorating  the morphisms of $\mbp_{{\mathbf A}_0}$. Let $(G, K,\alpha,\tau)$ be the   crossed module corresponding to $\mbg_1$.  For the decorated version of $\mbp_{{\mathbf A}_0}$ we take as objects just the objects of $\mbp$, and as morphisms
$$({\tilde\gamma}, h)\in  \Mor(\mbp_{{\mathbf A}_0})\times K  ,$$
where $K=\ker s_1$, with $s_1$ being the source map on the morphisms of $\mbg_2$. We define source and targets by
\begin{equation}\label{E:stdecor}
s({\tilde\gamma}, h)=s({\tilde\gamma})\quad\hbox{and}\quad t({\tilde\gamma}, h)=t({\tilde\gamma})t_1(h)^{-1},\end{equation}
composition by
 \begin{equation}\label{E:composedecor}
({\tilde\gamma}_2, h_2)\circ ({\tilde\gamma}_1, h_1)=\left({\tilde\gamma}_21_{\tau(h_1)}\circ {\tilde\gamma}_1, h_2h_1\right).\end{equation}
Let
\begin{equation}\label{E:defmbpabsdec}
\mbp_{{\mathbf A}_0}^{\rm dec} 
\end{equation}
be the category thus defined.
Next we define
an action of $\mbg_1$ on $\mbp_{{\mathbf A}_0}^{\rm dec}$ by
\begin{equation}\label{E:actiondecor}
({\tilde\gamma}, h)\cdot (h_1, g_1)=\left({\tilde\gamma} 1_{g_1}, \alpha(g_1^{-1})(h_1^{-1}h)\right),\end{equation}
for all $(h_1, g_1)\in K \times G$,
or, equivalently,
\begin{equation}\label{E:actiondecor2}
({\tilde\gamma}, h)\cdot\phi =\left({\tilde\gamma} 1_{s(\phi)}, \phi^{-1}h 1_{s(\phi)}\right).\end{equation}
Then 
\begin{equation}\label{E:defmbpabsdecbun}
\pi:\mbp_{{\mathbf A}_0}^{\rm dec}\to\mbbb\end{equation}
 is a principal categorical bundle with structure group $\mbg_1$. The proof is the same as for Theorem \ref{T:pathbundlecat}. Thus we have constructed a `decorated' version of a given principal categorical bundle with  a given categorical connection.

Let ${\tilde\gamma}_u$ denote, as usual, the ${\mathbf A}_0$-horizontal lift of $\gamma\in\Mor(\mbbb)$ with source $u$. Now suppose $k^*$ is a map from $\Mor(\mbp_{{\mathbf A}_1})$ to $K$ satisfying:
\begin{equation}\label{E:hstarconds}
\begin{split}
k^*({\tilde\gamma}_u1_{g_1}) &= \alpha(g_1^{-1} )k^*({\tilde\gamma}_u) \\
k^*\left({\tilde\delta}_v\circ{\tilde\gamma}_u\right)&=k^*\left({\tilde\gamma}_u\right)k^*\left({\tilde\delta}_v\right)
\end{split}
\end{equation}
for all ${\mathbf A}_0$-horizontal lift  ${\tilde\gamma}_u, {\tilde\delta}_v$ for which ${\tilde\gamma}_v\circ{\tilde\delta}_u$ is defined and for all $g_1\in G$.  (The notation $k^*$ is not meant to suggest any type of `pullback'.) The first equality corresponds to (\ref{E:hug1}).
 The second equality above corresponds to the equality (\ref{E:huhom}) noted earlier. 
 
 Define the horizontal lift of $\gamma\in \mbbb$ with initial point $u$ on the fiber over $s(\gamma)$ to be
\begin{equation}\label{E:defhorgen}
h^{\rm hor}_u(\gamma)= \left({\tilde\gamma}_u, k^*({\tilde\gamma}_u)\right),\end{equation}
where ${\tilde\gamma}_u$ is the ${\mathbf A}_0$-horizontal lift of $\gamma$ with source $u$.

We summarize  this construction in:

\begin{theorem}\label{T:absdecconn} Let $ {\mathbf A}_0$ be a connection on   a principal  categorical bundle $\pi:\mbp\to\mbbb$ with structure categorical group $\mbg_0$, and let 
$${\mathbf P}_{{\mathbf A}_0}\to {\mathbf B}$$
be the corresponding categorical bundle where the morphisms  of  ${\mathbf P}_{{\mathbf A}_0}$ are the ${\mathbf A}_0$-horizontal morphisms of ${\mathbf P}$. Let $(G, K,\alpha,\tau)$ be the   crossed module corresponding to a categorical group $\mbg_1$, where the object group of $\mbg_1$ is the same as the object group $G$ of $\mbg_0$.  Then the construction described for (\ref{E:defmbpabsdecbun}) produces a categorical principal bundle  
\begin{equation}\label{E:PGB12}
\pi:\mbp_{{\mathbf A}_0}^{\rm dec}\to\mbbb,
\end{equation}
 with structure group ${\mathbf G}_1$. 

Let $k^*:\Mor(\mbp_{{\mathbf A}_0})\to K$ satisfy the conditions (\ref{E:hstarconds}). Then $h^{\rm hor}:\Mor(\mbbb)\to \Mor(\mbp_{{\mathbf A}_0}) $ defined by (\ref{E:defhorgen}) specifies a categorical connection  ${\mathbf A}_1$ on the categorical principal bundle  (\ref{E:PGB12}).
\end{theorem}
We omit the proof, which is a straightforward  abstract reformulation of  the first part of the proof of Theorem \ref{T:connectpapbt}.

\section{Path categories}\label{S:pathcat}

In this section we set up some definitions and conventions for categories whose morphisms are paths in some spaces. For the sake of avoiding technicalities that obscure the main ideas we will focus only on paths that are $C^\infty$ and constant near the endpoints. The latter condition makes it possible to compose two paths without losing smoothness at the interface of composition between the paths.

Our first focus is on a generalization of   pathspaces. Let $X$ be a smooth manifold.  Consider a `box'  
$$I=\prod_{k=1}^N[a_k,b_k]\subset\mbr^N,$$
where $N$ is some positive integer, and $a_k,b_k\in\mbr$ with $a_k<b_k$.  We denote by
$$C^\infty_c(I;X)$$
the set of all $C^\infty$ maps $I\to X$ with the property that  there is some $\epsilon>0$ such that for any $u_i\in [a_i,b_i]$ for $i\in\{1,\ldots,N\}-\{k\}$, the function $u_k\mapsto f(u_1,\ldots, u_N)$   is constant when $|u_k-a_k|<\epsilon$ and when $|u_k-b_k|<\epsilon$.   

Our interest in the following discussion is in the cases $N=1$ and $N=2$. 

Let $I_0$ be the `lower face' of $I$, and $I_1$ the `upper face' :
\begin{equation}\label{E:I0I1}
I_0=\left(\prod_{k=1}^{N-1}[a_k,b_k]\right)\times \{a_N\}\quad\hbox{and}\quad I_1= \left(\prod_{k=1}^{N-1}[a_k,b_k]\right)\times \{b_N\}.
\end{equation}

Note that if $f\in C^\infty_c(I;X)$ then  the restrictions $f|I_0$ and $f|I_1$ specify elements  in $C_c^\infty(\prod_{k=1}^{N-1}[a_k,b_k];X)$, and we think of these elements as the `initial'  (source) and `final' (target)  values of $f$.
Thus, we think of  $f\in C_c^\infty(I;X)$ as a `morphism' from its source $s(f)$ to its target $t(f)$:
\begin{equation}\label{E:sftf01}
s(f)(x)=f(x,a_N)\quad\hbox{and}\quad t(f)(x)=f(x,b_N).\end{equation}

Now consider a second box $J\subset\mbr^N$, whose lower face $J_0$ is the upper face $I_1$ of $I$:
$$J_0=I_1.$$
In particular, the largest $N$-th coordinate for points in $I$ is equal to the smallest  $N$-th coordinate for $J$ and $I\cup J$ is a box in $\mbr^N$.
Then we compose $f\in   C_c^\infty(I;X)$ with $g\in C_c^\infty(J;X)$ if
$t(f)=s(g)$, defining the composite $g\circ_{\rm v} f$ to be the element of $C^\infty_c(I\cup J; X)$ given by
\begin{equation}\label{E:fcircgdef}
g\circ_{\rm v} f(u) =\begin{cases} f(u) & \hbox{if $u\in I$;}\\
g(u) &\hbox{if $v\in J$.}
\end{cases}
\end{equation}
Clearly, the composition operation is associative.

Since  one should be able to compose a morphism with another if the target of the first is the source of the second, regardless of the exact domains of the morphisms when taken as maps, 
\begin{quote}\label{E:identify}
{\em  we identify $f\in C_c^\infty(I;X)$ and  $h\in C_c^\infty(K;X)$ if there is some $d\in\mbr^N$ such that $I=K+d$ and  $f(u)=h(u+d)$ for all $u\in K$. }
\end{quote}
We denote the resulting quotient set by
\begin{equation}\label{E:MapNX}
{\rm Map}_N(X)
\end{equation}
if the domains of the functions are boxes in $\mbr^N$.  Then there are well-defined source and target maps
\begin{equation}\label{E:stMapNX}
s, t:{\rm Map}_N(X)\to{\rm Map}_{N-1}(X)\end{equation}
for all positive integers $N$.  Clearly the composition $g\circ_{\rm v} f$ is meaningful for $f, g\in {\rm Map}_N(X)$ if $t(f)=s(g)$.  
In order to get a category we must have identity morphisms and so we quotient one more time, by requiring that 
\begin{equation}\label{E:figigf}
f\circ_{\rm v} i=f\qquad\hbox{and}\qquad i\circ_{\rm v} g=g\end{equation}
whenever the compositions $f\circ_{\rm v} i$ and $i\circ_{\rm v} g$ are meaningful and  $i\in C^\infty_c(I;X)$ is independent of the $N$-th coordinate direction.  Let
\begin{equation}
{\rm Mor}_N(X)
\end{equation}
be the set of all equivalence classes of $f\in {\rm Map}_N(X)$, with the equivalence being defined by the requirements (\ref{E:figigf}).  Then the source and target maps $s$ and $t$ in (\ref{E:stMapNX}) descend to well-defined maps
\begin{equation}\label{E:stMorN0X}
s_{\rm v}, t_{\rm v}:{\rm Mor}_N(X)\to{\rm Mor}_{N-1}(X).\end{equation}
Again, composition $g\circ_{\rm v} f$ is meaningful whenever $t(f)=s(g)$, and the operation of composition is associative.  For any $a\in {\rm Mor}_{N-1}(X)$ (here $N\geq 1$) let  $1_a\in {\rm Mor}_N(X)$ be the element with $a$ as both source and target, and $1_a$, viewed as a mapping on an $N$-box, is constant along the $N$-th coordinate direction (thus corresponding to the $i$ in (\ref{E:figigf})). Then $f\circ_{\rm v} 1_a=f$ and $g=1_a\circ_{\rm v} g$ whenever $s(f)=a$ and $t(g)=a$.

Thus, in this way we obtain, for every positive integer $N$, a category
\begin{equation}\label{E:cat0X}
\mbp_N(X)
\end{equation}
whose object set is ${\rm Mor}_{N-1}(X)$ and whose morphism set is ${\rm Mor}_{N}(X)$.

Our main interest is in the notion of parallel-transport.  Parallel-transport has certain invariance properties: for example, parallel-transport is invariant under reparametrization of paths and backtrack-erasure.  One way of expressing these properties is to say that parallel-transport is well-defined on categories whose morphisms are obtained by identification of certain classes of morphisms in $\mbp_N(X)$, such as those that are `thin homotopy' equivalent (as in Theorem \ref{T:thinhomot}). However, instead of passing to such quotient categories we have and will state the relevant invariance properties as they arise. The only quotients/identification we work with are the bare minimum necessary ones to ensure that a category is formed by the maps of interest.

In the following   we use terminology introduced in the context of (\ref{E:I0I1}).  Let $X$ be a manifold, $N\geq 2$ a positive integer,
$I=\prod_{j=1}^N[a_j,b_j]$ a box in $\mbr^N$, and $f:I\to X$ a $C^\infty$ map.  Let $C$ be a $C^\infty$ $(N-1)$-form on $X$ with values in the Lie algebra $L(H)$ of a Lie group $H$.  Consider the solution
$$w_f:[a_N,b_N]\to H: u\mapsto w(u)$$
to the diferential   equation
\begin{equation}\label{E:defwmulti}
w_f(u)^{-1}w_f'(u) =-\int_{\prod_{j=1}^{N-1}[a_j,b_j]}C\bigl(\partial_1f(t,u),\ldots, \partial_{N-1}f(t,u)\bigr)\,dt, 
\end{equation}
 with initial condition 
 $$w(a_N)=e.$$
We define
\begin{equation}\label{E:wCdef}
w_C(f)=w_f(b_N).
\end{equation}
With this notation we have the following result on composites:

\begin{prop}\label{P:ptmulti} Let $X$ be a manifold, $N\geq 2$ a positive integer,
$I $ and $J $ boxes in $\mbr^N$ such that the lower face of $J$ is the upper face of $I$.  Let $f\in C^\infty_c(I;X)$ and $g\in C_c^\infty(J;X)$ be  such that the composite $g\circ_{\rm v}f$, given by (\ref{E:figigf}),  is defined.  Let $C$ be a $C^\infty$ $(N-1)$-form on $X$ with values in the Lie algebra $L(H)$ of a Lie group $H$. Then
\begin{equation}\label{E:fgcompC}
w_C(g\circ_{\rm v}f)=w_C(f)w_C(g). 
\end{equation}
\end{prop}

The property (\ref{E:fgcompC})  makes it possible to construct examples of categorical connections over pathspaces; we have used it in (\ref{E:huhom}) and will make use of it again in a more complex setting in section \ref{S:catconpath}.

\noindent\underline{Proof}.  Let $I$ be the box $\prod_{j=1}^N[a_j,b_j]$ and $J$ the box $\prod_{j=1}^N[c_j,d_j]$; since the upper face of $I$ is the lower face of $J$, we have $b_N=c_N$ and $[a_i,b_i]=[c_i,d_i]$ for $i\in\{1,\ldots, N-1\}$. The function of $u$  on the right in (\ref{E:defwmulti}) is $C^\infty$ and so the differential equation (\ref{E:defwmulti}) has a $C^\infty$ solution that is uniquely deteremined by the initial value $w_f(a_N)$. Furthermore, for any $h\in H$ the left-translate $hw_f$ is the solution  of (\ref{E:defwmulti}) with initial value $h$. Hence $w_f(b_N)w_g$ is the solution of  the differential equation
$$w(u)^{-1}w'(u) =-\int_{\prod_{j=1}^{N-1}[c_j,d_j]}C\bigl(\partial_1g(t,u),\ldots, \partial_{N-1}g(t,u)\bigr)\,dt, $$
with initial condition $w(c_N)=w_f(b_N)$. Thus the composite
\begin{equation}
\bigl(w_f(b_N)w_g\bigr)\circ w_f: [a_N, d_N]\to H: u\mapsto\begin{cases}  w_f(u) &\hbox{if $u\in [a_N,b_N]$;}\\
w_f(b_N)w_g(u) &\hbox{if $u\in [c_N,d_N]$}\end{cases}
\end{equation}
is continuous (with value $w_f(b_N)$ at $u=b_N$), has initial value $w_f(a_N)=e$ and is a solution of the differential equation 
\begin{equation}\label{E:wprucomp}
w(u)^{-1}w'(u) =-\int_{\prod_{i=1}^{N-1}[a_i,b_i]}C\bigl(\partial_1({g\circ_{\rm v} f})(t,u),\ldots, \partial_{N-1}({g\circ_{\rm v} f})(t,u)\bigr)\,dt, 
\end{equation}
at all $u\in [a_N,d_N]$ except possibly at $u=b_N$. Thus  $w_f(b_N)w_g$ agrees with the solution $w_{{g\circ_{\rm v} f}}$ of (\ref{E:wprucomp}), with initial value $e$, at all points $u\in [a_N,d_N]$ except possibly at $b_N$;   continuity of both   $w_f(b_N)w_g$ and  $w_{{g\circ_{\rm v} f}}$ ensures then that these are equal also at $u=b_N$. Thus
$$ w_{{g\circ_{\rm v} f}}=w_f(b_N)w_g(u)\qquad\hbox{for all $u\in [a_N,d_N]$.}$$
Taking $u=d_N$ gives
$$w_{{g\circ_{\rm v} f}}(d_N)=w_f(b_N)w_g(d_N),$$
which is just the identity (\ref{E:fgcompC}).
\fbox{QED}
 
 We can modify $w_C$ to another example of a function with a property similar to (\ref{E:fgcompC}):
 
 \begin{prop}\label{P:ptmulti2} Let $X$ be a manifold, $N\geq 2$ a positive integer.  Let $C$ be a $C^\infty$ $(N-1)$-form on $X$ with values in the Lie algebra $L(H)$ of a Lie group $H$, and  $w_0$  any $H$-valued function on $\Mor_{N-1}(X)$:
 \begin{equation}\label{E:defw0}
 w_0:   \Mor_{N-1}(X)\to H.
 \end{equation}
 Define
\begin{equation}\label{E:wCDf}
w_{C,0}(f)=w_0\bigl(s(f)\bigr)w_C(f)w_0\bigl(t(f)\bigr)^{-1}
\end{equation}
for all $f\in C^\infty_c(I;X)$, where $I$ is any box in $\mbr^N$. Then
\begin{equation}\label{E:fgcompC0}
w_{C,0}(g\circ_{\rm v}f)=w_{C,0}(f)w_{C,0}(g). 
\end{equation}
for all $f\in C^\infty_c(I;X)$ and $g\in C_c^\infty(J;X)$ for which  $g\circ_{\rm v}f$ is defined, with $I$ and $J$ being boxes in $\mbr^N$.
\end{prop}

Note that in particular we could take for $w_0$ the function $w_D$ obtained from an $(N-2)$-form $D$ on $X$ with values in $H$.

\noindent\underline{Proof}. The composite $g\circ_{\rm v}f$  has source $s(f)$ and target $t(g)$; hence
 \begin{equation}
\begin{split}
w_{C,0}(g\circ_{\rm v}f) &=w_0\bigl(s(g\circ_{\rm v}f)\bigr)w_C(g\circ_{\rm v}f) w_0\bigl(t(g\circ_{\rm v}f)\bigr)^{-1}\\
&\hskip -.25in = w_0\bigl(s(f)\bigr)w_C(f)w_C(g)w_0\bigl(t(g)\bigr)^{-1}\\
&\hskip -.25in= w_0\bigl(s(f)\bigr)w_C(f)w_D\bigl(t(f)\bigr)^{-1}\, w_0\bigl(t(f)\bigr) w_C(g)w_0\bigl(t(g)\bigr)^{-1}\\
&\qquad\qquad\hbox{(using $t(f)=s(g)$.)}\\
&\hskip -.25in= w_0\bigl(s(f)\bigr)w_C(f)w_0\bigl(t(f)\bigr)^{-1}\, w_0\bigl(s(g)\bigr) w_C(g)w_0\bigl(t(g)\bigr)^{-1}\\
&\hskip -.25in= w_{C,0}(f)w_{C,0}(g).
\end{split}
\end{equation}

\section{Categorical connections over pathspace}\label{S:catconpath}

We will construct a categorical connection over a pathspace using the $1$-form $\oab$ given in (\ref{E:defoab2}).  For this purpose it will be more convenient to work with $C^\infty$ paths that are constant near their endpoints, and define paths of paths as in section \ref{S:pathcat}.  As before we work with a Lie crossed module  $(G, H, \alpha,\tau)$, connection forms $A$ and $\ovA$ on a principal $G$-bundle $\pi:P\to M$, and an $L(H)$-valued $\alpha$-equivariant $2$-form $B$ on $P$ with values in $L(H)$ that vanishes on $(v,w)$ whenever $v$ or $w$ is a vertical vector in $P$.

Let  ${\mathbf P}_1(X)$ and ${\mathbf P}_2(X)$ be the categories described in section \ref{S:pathcat}. We will work with the case where $X$ is either $M$ or $P$.

In Theorem  \ref{T:connectpapbt} we constructed a categorical connection on a categorical bundle whose objects are points and whose morphisms are paths.  We now construct an example that is an analog of this, but one dimension higher in the sense that the objects are now paths and the morphisms are paths of paths.  We focus on the subcategory  ${\mathbf P}^{\ovA}_1(P)$ of  ${\mathbf P}_1(P)$ in which  the morphisms arise from $\ovA$-horizontal paths on $P$, and the subcategory   ${\mathbf P}^{\oab}_2(P)$ of ${\mathbf P}_2(P)$ where the morphisms come from $\oab$-horizontal paths of $\ovA$-horizontal paths on $P$.  We have then the categorical principal bundle

\begin{equation}\label{E:P02PM}
\pi:{\mathbf P}^{\oab}_2(P)\to {\mathbf P}_2(M),\end{equation}
whose structure categorical group is the categorical group ${\mathbf G}_d$  (object set is $G$ and morphisms are all the identity morphisms). This is the analog of Example P3 in section \ref{S:categorbundle}, one dimension higher, with paths replaced by paths of paths. The action of   any object $g\in G$ on any object ${\tilde\gamma}$ of ${\mathbf P}^{\oab}_2(P)$ produces ${\tilde\gamma}g$, which is again an $\ovA$-horizontal path on $P$. The action of $1_g:g\to g$ on ${\tilde\Gamma}\in \Mor\left({\mathbf P}^{\oab}_2(P)\right)$ produces ${\tilde\Gamma}g$:
$${\tilde\Gamma}1_g\stackrel{\rm def}{=}{\tilde\Gamma}g,$$
  which arises from the map $[s_0,s_1]\times[t_0,t_1]\to P:(s,t)\mapsto {\tilde\Gamma}(s,t)g$, where we use  a representative map ${\tilde\Gamma}$.  Note that if ${\tilde\Gamma}_s={\tilde\Gamma}_{s_0}$ for all $s\in [s_0,s_1]$ then ${\tilde\Gamma}_sg={\tilde\Gamma}_{s_0}g$ for all $s\in [s_0,s_1]$.

In order to produce a categorical connection on the bundle (\ref{E:P02PM}) we have to provide the rule for horizontal lifts of morphisms in ${\mathbf P}_2(M)$ to morphisms in $ {\mathbf P}^{\oab}_2(P)$. A morphism in ${\mathbf P}_2(M)$ arises from some $\Gamma:[s_0,s_1]\times [t_0,t_1]\to M:(s,t)\mapsto\Gamma_s(t)$; let ${\tilde\Gamma}^{\rm h}:[s_0,s_1]\times[t_0,t_1]\to P:(s,t)\mapsto {\tilde\Gamma}^{\rm h}_s(t)$
be its $\oab$-horizontal lift, with a given choice of initial path ${\tilde\Gamma}^{\rm h}_{s_0}$. It is readily checked that there is an $\epsilon>0$ such that each path ${\tilde\Gamma}^{\rm h}_s:[t_0,t_1]\to P$ is   constant within distance $\epsilon$ from $t_0$ and from $t_1$  (because the same is true for the projected path ${\Gamma}_s$ of which ${\tilde\Gamma}^{\rm h}_s$ is an $\ovA$-horizontal lift). Moreover, ${\tilde\Gamma}^{\rm h}_s={\tilde\Gamma}^{\rm h}_{s_0}$ for $s$ near $s_0$, and ${\tilde\Gamma}^{\rm h}_s={\tilde\Gamma}^{\rm h}_{s_1}$ for $s$ near $s_1$. Furthermore, of course, ${\tilde\Gamma}^{\rm h}$ is $C^\infty$. Hence ${\tilde\Gamma}^{\rm h}\in \Mor\left({\mathbf P}^{\oab}_2(P)\right)$.  We label this with the given initial path ${\tilde\gamma}={\tilde\Gamma}^{\rm h}_{s_0}$:
$${\tilde\Gamma}^{\rm h}_{{\tilde\gamma}}.$$
It is   clear that the assignment
$$\Gamma\mapsto {\tilde\Gamma}^{\rm h}_{{\tilde\gamma}}$$
is functorial in the sense that composites of paths are lifted to composites of the lifted paths; this is the analog of (\ref{E:horliftfunct}), and  says that the horizontal lift of ${\Delta}\circ_{\rm v}\Gamma$ with initial path ${\tilde\gamma}$ is
\begin{equation}\label{E:horliftfunct2}
\left({\widetilde{{\Delta}\circ_{\rm v}\Gamma}}\right)^{\rm h}_{\tilde\gamma} = {\tilde \Delta}^{\rm h}_{{\tilde\delta}}\circ_{\rm v} {\tilde \Gamma}^{\rm h}_{{\tilde\gamma}},\end{equation}
where
$${\tilde\delta}=s\left({\tilde\Delta}^{\rm h}_{{\tilde\delta}}\right)=t\left({\tilde\Gamma}^{\rm h}_{{\tilde\gamma}}\right).$$
(see (\ref{E:fcircgdef}) for the definition of `vertical' composition $\circ_{\rm v}$).
The condition of `rigid motion' (\ref{E:horliftrightact}) is also clearly valid since $s\mapsto {\tilde\Gamma}_sg$ is $\oab$-horizontal if $s\mapsto {\tilde\Gamma}_s$ is $\oab$-horizontal (by (\ref{E:oABvg})). $\square$

    \section{Parallel-transport for decorated paths}\label{S:ptdp}

       We turn now to our final example of a categorical connection. This will provide a connection on a decorated bundle over space of paths.  Before turning to the technical details let us summarize the essence of the construction. As input we have a connection $\ovA$ on a principal $G$-bundle
$$\pi:P\to M.$$
Next let $\mbg_1$ be a categorical group with associated crossed module
$$(G, H, \alpha_1,\tau_1).$$
 Using a connection $A$ and an equivariant $2$-form $B$ on $P$ with values in $L(H)$,  we have the $1$-form
 $$\oab$$ 
 as specified earlier in (\ref{E:defoab2}).  Now consider a categorical group $\mbg_2$ with
 $$\Obj(\mbg_2)=\Mor(\mbg_1),$$
 with associated crossed module
 $$(H\rtimes_{\alpha_1}G, K, \alpha_2,\tau_2).$$
We will use this to construct a doubly-decorated category $\mbp_2^{\rm dec}(P)^{\rm dec}_{\oab}$ whose objects are of the form
 $$({\tilde\gamma}, h),$$
 where $\tilde\gamma$ is any $\ovA$-horizontal path on  $P$ and $h\in H$, and whose morphisms are of the form
 $$({\tilde\Gamma}, h, k),$$
 where $\tilde\Gamma$ is any $\oab$-horizontal path of paths on $P$.  We will show that there is a principal categorical bundle
 $$\mbp_2^{\rm dec}(P)^{\rm dec}_{\oab}\to\mbp_2(M),$$
 with  structure categorical group $\mbg_2$, and then construct categorical connections on this bundle. Briefly put, our method provides a way of constructing parallel-transport of decorated paths
 $$({\tilde\gamma}, h)$$
 along doubly decorated  paths of paths
 $$({\tilde\Gamma}, h, k),$$
 where $\tilde\gamma$ is any $\ovA$-horizontal path on the original bundle and ${\tilde\Gamma}$ is an $\oab$-horizontal path on the bundle of $\ovA$-horizontal paths over the pathspace of $M$.

 Let us recall the construction provided by Theorem   \ref{T:absdecconn}.  Consider a categorical connection ${\mathbf A}_0$ on a categorical principal bundle ${\mathbf P}\to {\mathbf B}$, with structure categorical group ${\mathbf G}_0$. Let ${\mathbf P}_{{\mathbf A}_0 }\to {\mathbf B}$ be the categorical bundle obtained by working only with ${\mathbf A}_0$-horizontal morphisms of  ${\mathbf P}$.  Next consider a map
 $$k^*: \Mor\Bigl({\mathbf P}_{{\mathbf A}_0 }\Bigr)\to K,$$
 where $K$ is a group,
 satisfying  
 \begin{equation}\label{E:kstarconds}
\begin{split}
k^*({\tilde\gamma}_u1_{g_0}) &= \alpha({g_0}^{-1} )k^*({\tilde\gamma}_u) \\
k^*\left({\tilde\delta}_v\circ{\tilde\gamma}_u\right)&=k^*\left({\tilde\gamma}_u\right)k^*\left({\tilde\delta}_v\right)
\end{split}
\end{equation}
for all $g_0\in  \Obj({\mathbf G}_0)$ and  all ${\tilde\delta}_v, {\tilde\gamma}_u \in \Mor\Bigl({\mathbf P}_{{\mathbf A}_0 }\Bigr)$ that are composable.   Now let
$${\mathbf G}_1$$
be a categorical group whose object group is the same as the object group of $\mbg_0$:
\begin{equation}
\Obj(\mbg_1)=\Obj(\mbg_0)=G.
\end{equation}
Theorem   \ref{T:absdecconn} then provides   a `decorated' categorical principal bundle ${\mathbf P}_{{\mathbf A}_0}^{\rm dec}\to {\mathbf B}$, with structure group ${\mathbf G}_1$ (whose objects are the objects of ${\mathbf G}_0$) along with a categorical connection ${\mathbf A}_1$ on this bundle. Thus, the objects of  ${\mathbf P}_{{\mathbf A}_0}^{\rm dec}$ are the objects ${\tilde\gamma}$ of ${\mathbf P}_{{\mathbf A}_0}$, but the morphisms are decorated morphisms of $\mbp_{{\mathbf A}_0}$:
$$({\tilde\gamma}, h)\in \Mor(\mbp_{{\mathbf A}_0})\times H.$$

We now apply this procedure   with $\mbp \to \mbbb$ being a categorical bundle 
$${\mathbf P}^{\rm dec}_2(P)\to {\mathbf P}_2(M),$$
which we now describe. 

Let $\pi:P\to M$  a principal $G$-bundle equipped with connection $\ovA$.  Then we have a categorical principal bundle
$$\mbp^{\ovA}_1(P)\to \mbp_1(M),$$
where the object set  of $\mbp^{\ovA}_1(P)$ is $P$ and the morphisms arise from $\ovA$-horizontal paths ${\tilde\gamma}$ on $P$. The structure categorical group $\mbg_0$  is discrete, with object group being $G$.  There is a categorical connection $\mbba_0$ on this bundle: it associates to any $\gamma\in \Mor\bigl( \mbp_1(M)\bigr)$ the $\ovA$-horizontal lift ${\tilde\gamma}_u$ through any given initial point $u\in P$ in the fiber over $s(\gamma)$.

Now let  ${\mathbf G}_1$ be a categorical Lie group with crossed module $(G, H,\alpha_1,\tau_1)$. Then by  Theorem  \ref{T:absdecconn} we have  the decorated construction, yielding a categorical principal bundle
$${\mathbf P}_1^{\ovA}(P)^{\rm dec}\to {\mathbf P}_1(M),$$
whose object set is $P$ and whose morphisms  are of the form
$$({\tilde\gamma}, h),$$
where ${\tilde\gamma}$ arises from an $\ovA$-horizontal path on $P$, and $h\in H$. The structure categorical group is ${\mathbf G}_1$.

 We define the category $ {\mathbf P}_2^{\rm dec}(P)_{\oab}$ as follows. Its objects are the morphisms $({\tilde\gamma}, h)$ of ${\mathbf P}_1^{\ovA}(P)^{\rm dec}$. A morphism of $ {\mathbf P}_2^{\rm dec}(P)_{\oab}$ is of the form
$$({\tilde\Gamma}, h)\in \Mor\bigl( {\mathbf P}_2^{\rm dec}(P)_{\oab}\bigr)\times H$$
where ${\tilde\Gamma}$ is $\oab$-horizontal, with source and target being
\begin{equation}\label{E:stGamtildh}
s({\tilde\Gamma}, h)=\bigl(s({\tilde\Gamma}), h\bigr),\quad\hbox{and}\quad  
t({\tilde\Gamma}, h)=\bigl(t({\tilde\Gamma}), h\bigr),
\end{equation}
and composition being
\begin{equation}\label{E:Gamhhcomp}
({\tilde\Gamma}_2,h)\circ ({\tilde\Gamma}_1,h)=({\tilde\Gamma}_2\circ_{\rm v}{\tilde\Gamma}_1, h).
\end{equation} 
Then 
$$ {\mathbf P}_2^{\rm dec}(P)_{\oab} \to {\mathbf P}_2(M) $$
is a  categorical principal bundle with structure categorical group $\mbg_0$ being discrete, with object group $H\rtimes_{\alpha_1}G$:
$$\Obj(\mbg_0)=\Mor(\mbg_1)= H\rtimes_{\alpha_1}G.$$
(We use the notation $\mbg_0$ again in order to make the application of Theorem \ref{T:absdecconn} clearer.)  The action of an object $(h_1,g_1)\in \Obj(\mbg_0)$ on an object $({\tilde\gamma},h)$ of  $ {\mathbf P}_2^{\rm dec}(P)_{\oab}$ is given by
\begin{equation}\label{E:doubdecaction}
({\tilde\gamma},h)\cdot (h_1, g_1) = \Bigl({\tilde\gamma}g_1, \alpha(g_1^{-1})(h_1^{-1}h)\Bigr).
\end{equation}

The action of the identity morphism $1_{(h_1,g_1)}\in\Mor(\mbg_0)$ on $({\tilde\Gamma},h)$ is given by
\begin{equation}\label{E:TidGamhh1g1ra}
({\tilde\Gamma},h)1_{(h_1,g_1)}=\Bigl({\tilde\Gamma}g_1, \alpha(g_1^{-1})(h_1^{-1}h)\Bigr).
\end{equation}
For the categorical connection ${\mathbf A}_0$ we take the lifting of $\Gamma$ through $({\tilde\gamma}, h)$ to be
$$({\tilde\Gamma},h)$$
where $\tilde\Gamma$ is   $\oab$-horizontal  and has initial path $s({\tilde\Gamma})={\tilde\gamma}$.   

Now let ${\mathbf G}_2$ be a categorical group for which
\begin{equation}
\Obj({\mathbf G}_2)=\Obj(\mbg_0)=\Mor(\mbg_1)=H\rtimes_{\alpha_1}G,
\end{equation}
with Lie crossed module
\begin{equation}
(H\rtimes_{\alpha_1}G, K,\alpha_2,\tau_2).
\end{equation}

Then by Theorem  \ref{T:absdecconn} we obtain a doubly decorated categorical bundle
\begin{equation}\label{E:P2Mdecom}
{\mathbf P}_2^{\rm dec}(P)^{\rm dec}_{\oab}\to {\mathbf P}_2(M),
\end{equation}
with structure categorical group ${\mathbf G}_2$, whose object group is $H\rtimes_{\alpha_1} G=\Obj({\mathbf G}_0)$ and whose morphism group is $K\rtimes_{\alpha_2}  (H\rtimes_{\alpha_1}  G)$.

A morphism of  ${\mathbf P}_2^{\rm dec}(P)^{\rm dec}_{\oab}$ is of the form
$$({\tilde\Gamma}, h, k),$$
where $\tilde\Gamma$ is $\oab$-horizontal, $h\in H$ and $k\in K$; its source and target are (obtained from (\ref{E:compmor})):
\begin{equation}
s({\tilde\Gamma}, h, k)=\bigl(s({\tilde\Gamma}), h)\qquad\hbox{and}\qquad 
t({\tilde\Gamma}, h, k)=\bigl(  t({\tilde\Gamma}), h\bigr)\tau_2(k^{-1}),
\end{equation}
where in the second term on the right note that $\tau_2(k^{-1})\in H\times_{\alpha_1}G$ acts on the right on $\Mor(\mbp_1^{\rm dec}(P))$.

The composition of morphisms in this decorated bundle is given (again from (\ref{E:compmor})) by
\begin{equation}\label{E:compdecmor}
({\tilde\Delta},h, k) \circ_{\rm v}\bigl({\tilde\Gamma},h', k') =\bigl(({\tilde\Delta}, h)\tau_2(k')\circ_{\rm v} ({\tilde\Gamma}, h'), kk'\bigr)
\end{equation}

The right action of $\mbg_2$ on ${\mathbf P}_2^{\rm dec}(P)^{\rm dec}_{\oab}$ is given as follows. On objects the action of $(h_1,g_1)\in H\rtimes_{\alpha_1}G$ on $({\tilde\gamma}, h)$ is
\begin{equation}\label{E:h1g1tidgamhra}
({\tilde\gamma}, h)(h_1,g_1)=\Bigl({\tilde\gamma}g_1, \alpha_1(g_1^{-1})\bigl(h_1^{-1}h\bigr)\Bigr),
\end{equation}
just as seen before in (\ref{E:actiondecor}). On morphisms, the right action of $(k_1, h_1, g_1)\in K\rtimes_{\alpha_2}(H\rtimes_{\alpha_1}G)$ on $({\tilde\Gamma}, h,k)$ is given by
\begin{equation}\label{E:h1g1tidGamhra}
({\tilde\Gamma}, h,k)(k_1,h_1,g_1)=\Bigl({\tilde\Gamma}g_1, \alpha_1(g_1^{-1})(h_1^{-1}h), \alpha_2\bigl((h_1g_1)^{-1}\bigr)(k_1^{-1}k)\Bigr)
\end{equation}

Now we will construct a categorical connection on the categorical principal bundle (\ref{E:P2Mdecom}). To this end assume that we are given a  map from the morphisms of the   bundle $ {\mathbf P}_2^{\rm dec}(P)_{\oab}$ (before the $K$-decoration) to the group $K$:
$$k^*: \Mor\Bigl( {\mathbf P}_2^{\rm dec}(P)_{\oab}\Bigr)\to K,$$
  satisfying the conditions (\ref{E:kstarconds}):
  \begin{equation}\label{E:kstarconds2}
  \begin{split}
  k^*({\tilde\Gamma}1_{g_1}) &=\alpha_2(g_1^{-1})k^*({\tilde\Gamma})\\
  k^*({\tilde\Delta} \circ_{\rm v}{\tilde\Gamma}) &=k^*({\tilde\Gamma})k^*({\tilde\Delta}) \\
  \end{split}
  \end{equation}
  for all $g_1\in G$, and all $\oab$-horizontal ${\tilde\Gamma}$ and ${\tilde\Delta}$ for which the composite ${\tilde\Delta} \circ_{\rm v}{\tilde\Gamma}$ is defined. We have seen in Proposition \ref{P:ptmulti}   how such $k^*$ may be constructed by using `doubly path-ordered' integrals of forms with values in the Lie algebra of $K$ as well as, in Proposition \ref{P:ptmulti2}, how to obtain additional examples by including a `boundary term' to such integrals.  In more detail, suppose $C_1$ is an $L(K)$-valued $1$-form on $P$, and $C_2$ an  $L(K)$-valued $2$-form on $P$ such that they are both $0$ when contracted on a vertical vector and are $\alpha_2$-equivariant in the sense that
  \begin{equation}
  \begin{split}
  C_1((R_g)_*v_1) &=\alpha_2(g)^{-1}C_1(v_1)\\
  C_2\bigl((R_g)_*v_1, (R_g)_*v_2\bigr) &=\alpha_2(g^{-1})C_2(v_1,v_2)
  \end{split}
  \end{equation}
  for all $g\in G$, $v_1, v_2\in T_pP$ with $p$ running over $P$.  For ${\tilde\Gamma}:[s_0,s_1]\times [t_0,t_1]\to P$ we define $w_{C_2}^*({\tilde\Gamma})$ to be $w(t_1)$, where $[t_0,t_1]\to G: t\mapsto w(t)$ solves
  $$w(t)^{-1}w'(t)=-\int_{s_0}^{s_1}C_2\Bigl(\partial_u{\tilde\Gamma}(u,t), \partial_v{\tilde\Gamma}(u,t)\Bigr)\,du$$
  with $w(t_0)=e\in K$. (See equation (\ref{E:defwmulti}).)  By   $\alpha_2$-equivariance of $C_2$ we then have
  \begin{equation}
  k_2^*({\tilde\Gamma}g_1)=\alpha_2(g_1^{-1})k_2^*({\tilde\Gamma})
  \end{equation}
  for all $g_1\in G$. Moreover, $k_2^*$ satisfies the second property of $k^*$ in (\ref{E:kstarconds2}) by Proposition \ref{P:ptmulti}. Now define $w_0({\tilde\gamma})$, for any path ${\tilde\gamma}:[s_0,s_1]\to P$, to be  $w_1(s_1)$, where $w_1$ solves
  $$w_1(u)^{-1}w_1'(u)= -C_1\bigl({\tilde\gamma}'(u)\bigr),$$
  with $w_1(s_0)=e\in K$. Then by equivariance of $C_1$ we have
  $$ w_1({\tilde\gamma} {g_1})  =\alpha_2(g_1^{-1})w_1({\tilde\gamma})$$
  for all $g_1\in G$. Finally, set
  \begin{equation}\label{E:defkstartwopart}
  k^*({\tilde\Gamma})=w_0\bigl(s({\tilde\Gamma})\bigr)w_{C_2}({\tilde\Gamma})
  w_0\bigl(t({\tilde\Gamma})^{-1}.
  \end{equation}
  Then $k^*$ clearly satisfies the first condition in (\ref{E:kstarconds2}) because both $w_0$ and $w_{C_2}$ satisfy this condition; moreover, $k^*$ also satisfies the second condition in (\ref{E:kstarconds2})  by Proposition \ref{P:ptmulti2}.

  To construct a connection on the decorated bundle we need to obtain a mapping $\kappa^*$  similar to $k^*$ but for the decorated bundle:
 $$\kappa^*: \Mor\Bigl({\mathbf P}_2^{\rm dec}(P)_{\oab}^{\rm dec}\Bigr)\to K,$$
  satisfying the conditions (\ref{E:kstarconds}):
  \begin{equation}\label{E:kapstarconds2}
  \begin{split}
  \kappa^*\bigl(({\tilde\Gamma},h) 1_{(h_1,g_1)}\bigr) &=\alpha_2((h_1g_1)^{-1})\kappa^*\bigl({\tilde\Gamma},h \bigr)\\
  \kappa^*\bigl(({\tilde\Delta},h) \circ_{\rm v}\bigl({\tilde\Gamma},h)) &=\kappa^*\bigl( {\tilde\Gamma},h \bigr)\kappa^*\bigl( {\tilde\Delta},h \bigr),
  \end{split}
  \end{equation}
  where, in the second relation, note that for the composition on the left side to exist the $h$-component must be the same for ${\tilde\Delta}$ and ${\tilde\Gamma}$. 
  
  \begin{lemma}\label{L:kapstar}  Suppose 
  $$k^*: \Mor\Bigl({\mathbf P}_2^{\rm dec}(P)_{\oab}\Bigr)\to K,$$
  satisfies the conditions (\ref{E:kstarconds2}). Then the mapping 
   $$\kappa^*: \Mor\Bigl({\mathbf P}_2^{\rm dec}(P)_{\oab}^{\rm dec}\Bigr)\to K,$$
specified by
  \begin{equation}\label{E:defkapstar}
  \kappa^*({\tilde\Gamma},h)=\alpha_2(h)\bigl(k^*({\tilde\Gamma})\bigr)
  \end{equation}
  satisfies (\ref{E:kapstarconds2}).
  \end{lemma}
  \noindent\underline{Proof}. 
Using the formula for the right action for decorated bundles given in (\ref{E:morrightactpath}),  we have
\begin{equation}\label{E:morrightactpath3}
({\tilde\Gamma}, h)1_{(h_1,g_1)}=\bigl({\tilde\Gamma}g_1, \alpha_1(g_1^{-1})(h_1^{-1}h)\bigr).
  \end{equation} 
  Applying $\kappa^*$ we have
  \begin{equation}
  \begin{split}
  \kappa^*\bigl(({\tilde\Gamma}, h)1_{(h_1,g_1)}\bigr) &= \kappa^*\bigl({\tilde\Gamma}g_1, \alpha_1(g_1^{-1})(h_1^{-1}h)\bigr)\\
  &\hskip -.5in =\alpha_2\bigl(\alpha_1(g_1^{-1})(h_1^{-1}h)\bigr)\bigl(k^*({\tilde\Gamma}g_1) \bigr)\\
  &\hskip -.5in =\alpha_2\bigl(\alpha_1(g_1^{-1})(h_1^{-1}h)\bigr)\bigl(\alpha_2(g_1^{-1})k^*({\tilde\Gamma}) \bigr) \quad\hbox{(using (\ref{E:kstarconds2}))} \\
  &\hskip -.5in =\alpha_2(g_1^{-1}h_1^{-1}h)\bigl(k^*({\tilde\Gamma}) \bigr)   \quad\hbox{(using Lemma \ref{L:G2H2al2})}\\
  &\hskip -.5in =\alpha_2\bigl((h_1g_1)^{-1}\bigr) \bigl(\alpha_2(h)\bigl(k^*({\tilde\Gamma}) \bigr)  \bigr)\\
    &\hskip -.5in =\alpha_2((h_1g_1)^{-1})\kappa^*\bigl({\tilde\Gamma},h \bigr),
  \end{split}
  \end{equation}
  which establishes the first of the relations (\ref{E:kapstarconds2}).
  
  For the second relation we have
  \begin{equation}
  \begin{split}
   \kappa^*\bigl(({\tilde\Delta},h) \circ_{\rm v}\bigl({\tilde\Gamma},h)\bigr) &=
    \kappa^*\bigl( ({\tilde\Delta} \circ_{\rm v} {\tilde\Gamma}, h\Bigr) \qquad\hbox{(using (\ref{E:Gamhhcomp}))}\\
    &\hskip -.75in =\alpha_2(h)\bigl(k^*( {\tilde\Delta} \circ_{\rm v} {\tilde\Gamma})\bigr)\quad\hbox{(using (\ref{E:defkapstar}))}\\
    &\hskip -.75in =\alpha_2(h) \Bigl(k^*( {\tilde\Gamma})k^*( {\tilde\Delta}) \quad\hbox{(using (\ref{E:kstarconds2}))} \\
&\hskip -.75in =\alpha_2(h) \Bigl(k^*( {\tilde\Gamma}) \Bigr)\alpha_2(h) \Bigl(k^*( {\tilde\Delta})\Bigr)  \\
&\hskip -.75in =\kappa^*\bigl( {\tilde\Gamma},h \bigr)\kappa^*\bigl( {\tilde\Delta},h \bigr). \hskip 2in\fbox{QED}
  \end{split}
  \end{equation}
  
  Combining all of this we obtain a categorical connection $\mathbf A_2$ on the doubly decorated principal bundle ${\mathbf P}_2^{\rm dec}(P)^{\rm dec}_{\oab}$ with structure group $\mbg_2$ whose objects are the morphisms of an initially given categorical group $\mbg_1$, whose object group  $G$ in turn is the structure group of the original principal bundle $\pi:P\to M$.  Specifically, given a path of paths $\Gamma$ on $M$, and an initial $\ovA$-horizontal path ${\tilde\gamma}$ on $P$ lying above $s({\tilde\Gamma})$, and an element $h\in H$, the result of parallel-transport of $({\tilde\gamma}, h)$ by  the connection $\mathbf A_2$ along $\Gamma$ is the target 
  $$t({\tilde\Gamma}, h, k),$$
  where ${\tilde\Gamma}$ is $\oab$-horizontal with $s({\tilde\Gamma})={\tilde\gamma}$, and $k=\kappa^*({\tilde \Gamma}, h)$.     

\section{Associated bundles}\label{S:asocpt}

We define a categorical vector space to be a category $\mbv$  analogously to a categorical group. Both $\Obj(\mbv)$ and $\Mor(\mbv)$ are vector spaces, over some given field $\mbf$.  The addition operation
$$\mbv\times \mbv\to\mbv$$
is a functor, as is the scalar multiplication
$${\mathbf F} \times\mbv\to\mbv,$$
where ${\mathbf F}$ is the category with object set $\mbf$ and  morphisms just the identity morphisms for each $a\in\mbf$.

If $\mbk$ is a categorical group and $\mbv$ a categorical vector space then a {\em representation} $\rho$ of $\mbk$ on $\mbv$ is a functor
$$\rho:\mbk\times \mbv\to \mbv: (k,v)\mapsto \rho(k)v $$
such that $\rho(k):\Obj(\mbv)\to \Obj(\mbv)$, for each $k\in\Obj(\mbk)$
and $\rho(\phi):\Mor(\mbv)\to\Mor(\mbv) $, for all $\phi\in\Mor(\mbk)$, are linear,
and $\rho$ gives a representation of the group $\Obj(\mbk)$ on the vector space $\Obj(\mbv)$ as well as a representation of $\Mor(\mbk)$ on $\Mor(\mbv)$. 

Let $\pi:\mbp\to\mbbb$  be a principal categorical bundle with group $\mbk$, and $\rho$ a representation of $\mbk$ on a categorical vector space $\mbv$.  Then we construct a twisted product 
$$\mbp\times_\rho\mbv$$
An object of this category is an equivalence class $[p,v]$  of pairs $(p,v)$, with two such pairs $(p',v')$ and $(p,v)$ considered equivalent if  $p'=p\rho(k)$ and $v'=\rho(k^{-1})v$ for some $k\in\Obj(\mbk)$.  A morphism of $\mbp\times_\rho\mbv$ is an equivalence class $[F,f]$  of pairs $(F,f)$, with $(F,f)$ considered equivalent to $(F',f')$ if
$F'=F\rho(\phi)$ and $f'=\rho(\phi^{-1})f$ for some morphism $\phi$ of $\mbk$.
The target for $[F,f]$ is $[t(F), t(f)]$, and source is $[s(F), s(f)]$. There is a well-defined projection functor
$$\mbp\times_\rho\mbv\to\mbp$$
taking $[p,v]$ to $\pi_P(p)$ and $[F,f]$ to $\pi_P(F)$. We view this as the {\em associated vector bundle} for the given structures.

Given a connection on $\mbp$ we can construct parallel transport on $\mbp\times_\rho\mbv$ as follows. Consider a morphism $f:x\to y$ in $\mbbb$ and $(p,v)\in \Obj(\mbp)\times\Obj(\mbv)$ such that $\pi_P(p)=x$. Then we define the {\em parallel transport}  of $[p,v]$ along $f$ to be $[t(F), v]$, where $F$ is the horizontal lift of $f$ through $p$. 

We apply the abstract constructions above to the  specific principal categorical bundles we have studied before.

Let  $\overline{{\mathbf A}_0}$ be the categorical connection from Example CC1, associated with a connection $\ovA$ on a principal $G$-bundle $\pi:P\to B$.  A morphism of  $\mbp_0\times_{\rho}\mbv$ is an equivalence class of ordered pairs/triples
$$((\gamma; p,q); v, w),$$
where ${ \gamma}$ is a backtrack erased path on $M$ from $\pi(p)$ to $\pi(q)$, and $v, w\in V$. We think of $v$ as being located at the source of $ \gamma$ and $w$ being located at the target of $ \gamma$.  The equivalence relation between such triples is given by
\begin{equation}\label{E:gabvequiv}
((\gamma; p,q); v, w) \sim \bigl((\gamma; pg,qg); g^{-1}(v,w)\bigr).\end{equation}
From a connection $\ovA$  on $\pi:P\to B$ we have a connection $\overline{\mathbf A}_0$ on $\pi:\mbp_0\to\mbbb$, and then the corresponding parallel-transport of $[p,v]$  along $\gamma$ results in $[q, v]$, where $q$ is the target of the $\ovA$-horizontal lift of $\ovA$ initiating at $p$. This is exactly in accordance with the traditional notion of parallel transport for associated bundles.

 {\bf Acknowledgments.}   ANS   acknowledges support from a Mercator Guest Professorship at the University of Bonn from the Deutsche Forschungsgemeinschaft, and thanks Professor S. Albeverio for discussions and his kind hospitality.

\end{document}